\documentclass{amsproc}

\input xy
\xyoption{all}

\newtheorem{theorem}{Theorem}[section]
\newtheorem{proposition}[theorem]{Proposition}

\theoremstyle{definition}
\newtheorem{definition}[theorem]{Definition}
\newtheorem{example}[theorem]{Example}

\theoremstyle{remark}
\newtheorem{remark}[theorem]{Remark}

 \DeclareMathOperator{\comp}{\#}
 \DeclareMathOperator{\id}{id}

\begin{document}

\title{Thin fillers in the cubical nerves of omega-categories}

\author{Richard Steiner}

\address{Department of Mathematics\\University of Glasgow\\University Gardens\\
Glasgow\\Scotland G12 8QW}

\email{r.steiner@maths.gla.ac.uk}



\begin{abstract}
It is shown that the cubical nerve of a strict omega-category is a
sequence of sets with cubical face operations and distinguished
subclasses of thin elements satisfying certain thin filler
conditions. It is also shown that a sequence of this type is the
cubical nerve of a strict omega-category unique up to isomorphism;
the cubical nerve functor is therefore an equivalence of
categories. The sequences of sets involved are the analogues of
cubical T-complexes appropriate for strict omega-categories.
Degeneracies are not required in the definition of these
sequences, but can in fact be constructed as thin fillers. The
proof of the thin filler conditions uses chain complexes and chain
homotopies.
\end{abstract}

\maketitle

\section{Introduction} \label{S1}

This paper is concerned with the cubical nerves of strict
$\omega$-categories. It has been shown in~\cite{ABS} that these
cubical nerves are sequences of sets together with face maps,
degeneracies, connections and compositions, subject to various
identities, and that the functor taking an $\omega$-category to
its cubical nerve is an equivalence of categories. In this paper
we give a more conceptual characterisation of cubical nerves: a
cubical nerve is a sequence of sets with cubical face operations
and distinguished subclasses of thin elements such that certain
shells and boxes have unique thin fillers and such that one simple
extra condition is satisfied. Here a shell is a configuration like
the boundary of a cube, and a box is a configuration like the
boundary of a cube with one face removed.

There are similar characterisations of the cubical and simplicial
nerves of $\omega$-groupoids in terms of Dakin's
T-complexes~\cite{Dak} due to Ashley and Brown and Higgins
(\cite{Ash}, \cite{BHAlg}, \cite{BHCub}), and there is a similar
characterisation of the simplicial nerves of $\omega$-categories
due to Verity~\cite{V} (see also~\cite{StrWeak}). In the results
for $\omega$-groupoids every box or horn has a unique thin filler
(a horn is the simplicial analogue of a box) and there are
degeneracy operations, but there are no requirements on shells;
the resulting structure is called a cubical or simplicial
$T$-complex. In Verity's result also, certain horns have unique
thin fillers and there are degeneracies, but there are no
requirements on shells. Our result is like Verity's and is
different from the results on $\omega$-groupoids because we do not
require all boxes to have unique thin fillers. Our result differs
from all the previous results because we require thin fillers for
shells instead of requiring degeneracies. We need thin fillers for
shells in order to construct connections (compare the work of
Higgins in~\cite{Higg}), and it is more economical to construct
degeneracies from shells as well. This method may also be better
for extensions to weak $\omega$-categories, where one expects thin
fillers to exist but not necessarily to be unique. Degeneracies
correspond to identities, and in weak $\omega$-categories one does
not necessarily want unique identities, so one may also not want
unique degeneracies. Thin fillers for shells could be a suitable
alternative.

The main result of this paper, Theorem~\ref{2N}, is stated in
Section~\ref{S2}. The cubical nerve functor is described in
Section~\ref{S3}, with some calculations postponed to
Section~\ref{S6}, and the reverse functor is described in
Section~\ref{S4}. In Section~\ref{S5} the functors are shown to be
inverse equivalences. The proof of the thin filler conditions in
Sections \ref{S3} and~\ref{S6} is along the same lines as the
proof for the simplicial case given by Street in~\cite{StrFill},
but it is simplified by the use of chain complexes as
in~\cite{Ste}. Section~\ref{S7} gives an example showing how the
nerve of an $\omega$-category can differ from the nerve of an
$\omega$-groupoid.

\section{Statement of the main result} \label{S2}

Our result concerns sequences of sets with cubical face
operations, which we call precubical sets. We define precubical
sets in terms of the precubical category, which is an analogue of
the simplicial category but without degeneracies.

\begin{definition} \label{2A}
The \emph{precubical category} is the category with objects $[0]$,
$[1]$, \dots\ indexed by the natural numbers and with a generating
set of morphisms
$$\check\partial_1^-,\check\partial_1^+,\ldots,
 \check\partial_n^-,\check\partial_n^+\colon[n-1]\to[n]$$
subject to the relations
$$\check\partial_j^\beta\check\partial_i^\alpha
 =\check\partial_{i+1}^\alpha\check\partial_j^\beta
 \ \text{for $i\geq j$}.$$
A \emph{precubical set} is a contravariant functor from the
precubical category to sets; equivalently, a precubical set~$X$ is
a sequence of sets $X_0,X_1,\ldots\,$ together with \emph{face
operations}
$$\partial_1^-,\partial_1^+,\ldots,
 \partial_n^-,\partial_n^+\colon X_n\to X_{n-1}$$
such that
$$\partial_i^\alpha\partial_j^\beta
 =\partial_j^\beta\partial_{i+1}^\alpha
 \ \text{for $i\geq j$}.$$
If $X$~is a precubical set then the members of~$X_n$ are called
\emph{$n$-cubes}. An operation from $n$-cubes to $m$-cubes induced
by a morphism in the precubical category is called a
\emph{precubical operation}. A \emph{morphism of precubical sets}
from~$X$ to~$Y$ is a sequence of functions from~$X_n$ to~$Y_n$
commuting with the face operations.
\end{definition}

The precubical operations from $n$-cubes to $(n-1)$-cubes
correspond to the $2n$ faces of dimension $(n-1)$ in a standard
geometrical $n$-cube, and the relations between them correspond to
pairwise intersections. More generally, it is well-known that the
precubical operations from $n$-cubes to $m$-cubes correspond to
the $m$-dimensional faces of a geometrical $n$-cube, as follows.

\begin{proposition} \label{2B}
The precubical operations from $n$-cubes to $m$-cubes have unique
\emph{standard} decompositions
$$\partial_{i(1)}^{\alpha(1)}\ldots\partial_{i(n-m)}^{\alpha(n-m)}$$
with $1\leq i(1)<i(2)<\ldots<i(n-m)\leq n$.
\end{proposition}

In order to handle thin $n$-cubes we introduce the following
terminology.

\begin{definition} \label{2C}
A \emph{stratification} on a precubical set~$X$ is a sequence of
sets $t_1 X,t_2 X,\ldots\,$ such that $t_n X$ is a subset
of~$X_n$. A precubical set~$X$ with a stratification $t_n X$ is
called a \emph{stratified precubical set}, and the members of $t_n
X$ are called \emph{thin $n$-cubes}. A \emph{morphism of
stratified precubical sets} is a morphism of precubical sets
taking thin $n$-cubes to thin $n$-cubes.
\end{definition}

Note that $0$-cubes are never thin.

The cubical analogue of a horn is a configuration of the kind got
from the boundary of an $n$-cube by removing one face, and it is
called an $n$-box. We also need the configuration corresponding to
the complete boundary of an $n$-cube, which is called an
$n$-shell. The precise definitions are as follows.

\begin{definition} \label{2D}
Let $X$ be a precubical set. Then an \emph{$n$-shell}~$s$ in~$X$
is a collection of $(n-1)$-cubes~$s_i^\alpha$, indexed by the $2n$
face operations~$\partial_i^\alpha$ on $n$-cubes, such that
$$\partial_i^\alpha s_j^\beta
 =\partial_j^\beta s_{i+1}^\alpha
 \ \text{for $i\geq j$}.$$
A \emph{filler} for an $n$-shell~$s$ is an $n$-cube~$x$ such that
$\partial_i^\alpha x=s_i^\alpha$ for all~$\partial_i^\alpha$.
\end{definition}

\begin{definition} \label{2E}
Let $X$ be a precubical set and let $\partial_k^\gamma$ be a face
operation on $n$-cubes. Then an \emph{$n$-box}~$b$ in~$X$
opposite~$\partial_k^\gamma$ is a collection of
$(n-1)$-cubes~$b_i^\alpha$, indexed by the $(2n-1)$ face
operations~$\partial_i^\alpha$ on $n$-cubes other
than~$\partial_k^\gamma$, such that
$$\partial_i^\alpha b_j^\beta
 =\partial_j^\beta b_{i+1}^\alpha
 \ \text{for $i\geq j$}.$$
A \emph{filler} for an $n$-shell~$b$ is an $n$-cube~$x$ such that
$\partial_i^\alpha x=b_i^\alpha$ for
$\partial_i^\alpha\neq\partial_k^\gamma$.
\end{definition}

If $s$~is an $n$-shell and
$$\theta
 =\partial_{i(1)}^{\alpha(1)}\ldots\partial_{i(p)}^{\alpha(p)}$$
is a non-identity precubical operation on $n$-cubes, not
necessarily in standard form, then there is a well-defined
operation of~$\theta$ on~$s$ given by
$$\theta s
 =\partial_{i(1)}^{\alpha(1)}\ldots\partial_{i(p-1)}^{\alpha(p-1)}
 s_{i(p)}^{\alpha(p)},$$
and if $x$~is a filler for~$s$ then $\theta x=\theta s$. Similar
remarks apply to an $n$-box~$b$ opposite~$\partial_k^\gamma$,
except that the operation~$\partial_k^\gamma$ is not defined
on~$b$.

In the cubical nerve of an $\omega$-category we will require
certain boxes and shells to have unique thin fillers. We will now
describe what happens in dimension~$2$. We regard a $0$-cube as an
object of a category, and we regard a $1$-cube~$e$ as a morphism
from $\partial_1^- e$ to $\partial_1^+ e$. A $2$-cube~$x$ can then
be viewed as two composites
$$\partial_1^- x\circ\partial_2^+ x,
 \partial_2^- x\circ\partial_1^+ x\colon
 \partial_1^-\partial_2^- x\to\partial_1^+\partial_2^+ x,$$
together with some kind of higher morphism or $2$-morphism between
the two composites; see Figure~\ref{F1}. A $1$-cube is thin if it
is an identity morphism, and a $2$-cube is thin if its
$2$-morphism is an identity $2$-morphism. A $2$-shell is like a
$2$-cube with the $2$-morphism omitted, and a $2$-box is like a
$2$-shell with one of the edge morphisms omitted.

\begin{figure}
 $$\xymatrix{
 \partial_1^+\partial_1^-x=\partial_1^-\partial_2^+ x
 \ar[rrr]^{\partial_2^+ x}&&&
 \partial_1^+\partial_2^+ x=\partial_1^+\partial_1^+ x\\
 &
 \ar@{=>}[dr]\\
 &&\\
 \partial_1^-\partial_1^- x=\partial_1^-\partial_2^- x
 \ar@<2em>[uuu]^{\partial_1^- x}
 \ar[rrr]_{\partial_2^-x}&&&
 \partial_1^+\partial_2^- x=\partial_1^-\partial_1^+ x
 \ar@<-2em>[uuu]_{\partial_1^+ x}}$$
 \caption{A $2$-cube}
 \label{F1}
\end{figure}

It is now clear that a $2$-shell~$s$ has a unique thin filler if
it is commutative, that is $s_1^-\circ s_2^+=s_2^-\circ s_1^+$,
and that it has no thin filler otherwise. A $2$-box therefore has
a unique thin filler if it can be extended to a commutative
$2$-shell in a unique way. For example, a $2$-box~$b$
opposite~$\partial_1^-$ such that $b_2^+$~is an identity has a
unique thin filler~$x$, given by $\partial_1^- x=b_2^-\circ
b_1^+$; see Figure~\ref{F2}.

\begin{figure}
 $$\xymatrix{
 \ar[rr]^{b_2^+}_{=}&&\\
 \\
 \ar@{.>}[uu]^{b_2^-\circ b_1^+}
 \ar[rr]_{b_2^-}&&
 \ar[uu]_{b_1^+}}$$
 \caption{A $2$-box opposite $\partial_1^-$ with a unique thin filler}
 \label{F2}
\end{figure}

In general a box~$b$ opposite~$\partial_k^\gamma$ is to have a
unique thin filler if certain of its faces are thin, and we will
now explain which faces are involved. We start by considering a
class of \emph{extreme precubical operations}, where a precubical
operation is called extreme if its standard decomposition
$\partial_{i(1)}^{\alpha(1)}\ldots\partial_{i(p)}^{\alpha(p)}$ is
such that the signs $(-)^{i(r)-r}\alpha(r)$ are constant. Thus the
extreme precubical operations on $2$-cubes are
$$\id,\partial_1^-,\partial_1^+,\partial_2^-,\partial_2^+,
 \partial_1^-\partial_2^-,\partial_1^+\partial_2^+.$$
The non-extreme precubical operations on $2$-cubes are
$\partial_1^-\partial_2^+$ and $\partial_1^+\partial_2^-$, which
are in a sense between $\partial_1^-\partial_2^-$ and
$\partial_1^+\partial_2^+$ (see Figure~\ref{F1}). We will say that
an extreme precubical operation is
\emph{opposite~$\partial_k^\gamma$} if its standard decomposition
has a factor~$\partial_k^{-\gamma}$. The extreme precubical
operations opposite~$\partial_k^\gamma$ are indexed by the subsets
of $\{1,\ldots,n\}$ containing~$k$; for example the extreme
precubical operations on $5$-cubes opposite~$\partial_3^-$ are the
composites $\phi'\partial_3^+\phi''$ such that
$\phi'\in\{\id,\partial_1^-,\partial_1^+\partial_2^+,\partial_2^+\}$
and
$\phi''\in\{\id,\partial_4^+,\partial_4^+\partial_5^+,\partial_5^-\}$.
By omitting the factor~$\partial_k^{-\gamma}$ from the standard
decomposition for an extreme precubical operation~$\theta$
opposite~$\partial_k^\gamma$ we get a precubical operation which
in a sense joins~$\theta$ to~$\partial_k^\gamma$. We say that an
operation of this type is \emph{complementary
to~$\partial_k^\gamma$}. The precubical operations complementary
to~$\partial_k^\gamma$ are indexed by the subsets of
$\{1,\ldots,n\}$ not containing~$k$; for example the precubical
operations on $2$-cubes complementary to~$\partial_1^-$ are
$\partial_2^+$~and $\id$. It turns out that an $n$-box~$b$
opposite~$\partial_k^\gamma$ has a unique thin filler if $\theta
b$ is thin for every non-identity precubical operation on
$n$-cubes complementary to~$\partial_k^\gamma$. Boxes satisfying
this condition will be called \emph{admissible}; they correspond
to the admissible horns of~\cite{StrSimp}. For example, a
$2$-box~$b$ opposite~$\partial_1^-$ is admissible if and only if
$b_2^+$~is thin, in which case $b$~has a unique thin filler as in
Figure~\ref{F2}. We summarise these ideas concisely as follows.

\begin{definition} \label{2F}
Let $\partial_k^\gamma$ be a face operation on $n$-cubes and let
$\theta$ be a precubical operation on $n$-cubes. Then $\theta$~is
\emph{complementary to~$\partial_k^\gamma$} if the standard
decomposition of~$\theta$ has no factor
$\partial_k^-$~or~$\partial_k^+$ and if the insertion
of~$\partial_k^{-\gamma}$ produces a standard decomposition
$$\phi=\partial_{i(1)}^{\alpha(1)}\ldots\partial_{i(p)}^{\alpha(p)}$$
such that the sign $(-)^{i(r)-r}\alpha(r)$ is constant.
\end{definition}

\begin{definition} \label{2H}
In a stratified precubical set a box opposite~$\partial_k^\gamma$
is \emph{admissible} if $\theta b$ is thin for every non-identity
precubical operation~$\theta$ complementary
to~$\partial_k^\gamma$.
\end{definition}

\begin{remark} \label{2G}
A face operation~$\partial_l^\delta$ is complementary
to~$\partial_k^\gamma$ if and only if $k\neq l$ and
$(-)^k\gamma=(-)^l\delta$.
\end{remark}

\begin{example} \label{2J}
Let $b$ be an $n$-box opposite~$\partial_k^\gamma$ such that
$\theta b$ is thin whenever the standard decomposition of~$\theta$
does not contain $\partial_{k-1}^\gamma$, $\partial_k^-$,
$\partial_k^+$ or~$\partial_{k+1}^\gamma$. (In the case $k=1$ we
take it to be automatically true that the standard decomposition
does not have a factor~$\partial_{k-1}^\gamma$; in the case $k=n$
we take it to be automatically true that the standard
decomposition does not have a factor~$\partial_{k+1}^\gamma$.)
Then the standard decomposition of an operation~$\theta$ such that
$\theta b$ is not thin must have a factor $\partial_{k-1}^\gamma$,
$\partial_k^-$, $\partial_k^+$ or~$\partial_{k+1}^\gamma$, so it
cannot be complementary to~$\partial_k^\gamma$. Therefore $b$~is
admissible.
\end{example}

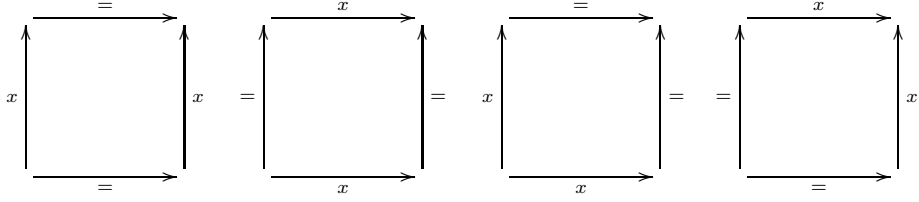
\begin{figure}
 $$\xymatrix{
 \ar[rr]^{=}&&&
 \ar[rr]^{x}&&&
 \ar[rr]^{=}&&&
 \ar[rr]^{x}&&\\
 \\
 \ar[uu]^{x}\ar[rr]_{=}&&\ar[uu]_{x}&
 \ar[uu]^{=}\ar[rr]_{x}&&\ar[uu]_{=}&
 \ar[uu]^{x}\ar[rr]_{x}&&\ar[uu]_{=}&
 \ar[uu]^{=}\ar[rr]_{=}&&\ar[uu]_{x}}$$
 \caption{Admissible $2$-shells}
 \label{F3}
\end{figure}

Next we consider thin fillers for shells. We have already observed
that a commutative $2$-shell has a unique thin filler, and a
similar result holds in general (see Proposition~\ref{3L}). In the
characterisation, however, we require a shell to have a unique
thin filler only when we can prove it to be commutative by using
admissible boxes. Shells of this kind will also be called
admissible. There are four types of admissible $2$-shell~$s$, as
shown in Figure~\ref{F3}. In each case there are two equal faces,
denoted~$x$ in the figure, which correspond to non-complementary
face operations, and omitting either of these faces produces an
admissible box. In general the definition is as follows.

\begin{definition} \label{2I}
In a stratified precubical set a shell~$s$ is \emph{admissible} if
there are distinct non-complementary face operations
$\partial_k^\gamma$~and~$\partial_l^\delta$ such that
$s_k^\gamma=s_l^\delta$ and such the boxes formed by removing
$s_k^\gamma$~or~$s_l^\delta$ from~$s$ are both admissible.
\end{definition}

\begin{example} \label{2K}
As in the first two cases of Figure~\ref{F3}, let $s$ be a shell
such that $s_k^-=s_k^+$ for some~$k$ and such that $\theta s$ is
thin whenever the standard decomposition of~$\theta$ does not
contain $\partial_k^-$~or~$\partial_k^+$. It follows from
Example~\ref{2J} that $s$~is admissible. Shells of this type will
be used to construct degeneracies.
\end{example}

\begin{example} \label{2L}
As in the last two cases of Figure~\ref{F3}, let $s$ be a shell
such that $s_k^\gamma=s_{k+1}^\gamma$ for some $k$~and~$\gamma$
and such that $\theta s$ is thin whenever the standard
decomposition of~$\theta$ does not contain $\partial_k^\gamma$,
$\partial_{k+1}^\gamma$ or
$\partial_k^{-\gamma}\partial_{k+1}^{-\gamma}$. Again it follows
from Example~\ref{2J} that $s$~is admissible. Shells of this type
will be used to construct connections.
\end{example}

For the admissible $2$-box~$b$ in Figure~\ref{F2}, note that the
additional face $b_2^-\circ b_1^+$ in the thin filler is thin if
all the $1$-cubes~$b_i^\alpha$ are thin. An analogous result holds
in general, and it is the final condition characterising nerves of
$\omega$-categories. We end up with the following structure.

\begin{definition} \label{2M}
A \emph{complete stratified precubical set} is a stratified
precubical set satisfying the following conditions: every
admissible box and every admissible shell has a unique thin
filler; if $x$~is the thin filler of an admissible box~$b$
opposite~$\partial_k^\gamma$ such that $b_i^\alpha$~is thin for
all $\partial_i^\alpha\neq\partial_k^\gamma$, then the additional
face $\partial_k^\gamma x$ is thin as well.
\end{definition}

The main result is now as follows.

\begin{theorem} \label{2N}
The cubical nerve functor is an equivalence between strict
$\omega$-categories and complete stratified precubical sets.
\end{theorem}

In the remainder of this paper we prove four results, Theorems
\ref{3N}, \ref{4M}, \ref{5D} and~\ref{5E}, whose conjunction is
equivalent to Theorem~\ref{2N}.

\begin{remark} \label{2O}
The only admissible boxes and shells used in the passage from
stratified precubical sets to the nerves of $\omega$-categories
are those of the types described in Examples \ref{2J}, \ref{2K}
and~\ref{2L}. It would therefore be sufficient to use these boxes
and shells in the definition of a complete stratified precubical
set.

Simplicial nerves behave in the same way. Indeed
Street~\cite{StrFill} shows that the simplicial nerve of an
$\omega$-category has a large class of admissible horns with
unique thin fillers, and Verity~\cite{V} shows that simplicial
nerves are characterised by a smaller class of complicial horns.
Admissible boxes in general correspond to admissible horns, and
the admissible boxes of Example~\ref{2J} correspond to complicial
horns.
\end{remark}

\begin{remark} \label{2P}
One of Dakin's axioms for a cubical or simplicial T-complex
(\cite{Ash}, \cite{BHCub}, \cite{Dak}) says that a thin $n$-cube
or $n$-simplex with all but one of its $(n-1)$-faces thin must
have its remaining $(n-1)$-face thin as well. In a complete
stratified precubical set we require this condition only when the
given thin faces form an admissible box. In the cubical nerve of
an $\omega$-category, it is in fact possible for a thin
$n$-cube~$x$ to have exactly one non-thin face~$\partial_i^\alpha
x$; it is even possible for this to happen when $x$~is the thin
filler of an admissible box, provided that this box is opposite an
operation other than~$\partial_i^\alpha$. An example is given in
Section~\ref{S7}. As before, one gets the same behaviour in
simplicial nerves.
\end{remark}

\section{From omega-categories to complete stratified precubical sets}
\label{S3}

In this section we show that the cubical nerve of an
$\omega$-category is a complete stratified precubical set. Recall
from \cite{ABS} that the cubical nerve of an $\omega$-category~$C$
is the sequence of sets
$$\hom(\nu I^0,C),\ \hom(\nu I^1, C),\ \ldots,$$
where $\nu I^n$ is an $\omega$-category associated to the
$n$-cube. We will establish the result by studying the
$\omega$-categories $\nu I^n$.

We begin by recalling the general theory of $\omega$-categories
\cite{StrSimp}. An $\omega$-category~$C$ is a set with a sequence
of compatible partially defined binary composition operations
subject to various axioms; in particular each composition
operation makes~$C$ into the set of morphisms for a category. We
will use $\comp_0$, $\comp_1$,~\dots\ to denote the composition
operations, and we will write $d_p^- x$ and $d_p^+ x$ for the left
and right identities of an element~$x$ under~$\comp_p$. We will
use the following properties.

\begin{proposition} \label{3A}
Let $C$ be an $\omega$-category. If $x\in C$ and $p<q$ then
$$d_p^\alpha d_q^\beta x=d_q^\beta d_p^\alpha x=d_p^\alpha x;$$
if $x\comp_q y$ is a composite and $p\neq q$ then
$$d_p^\alpha(x\comp_q y)=d_p^\alpha x\comp_q d_p^\alpha y;$$
the identities for~$\comp_p$ form a sub-$\omega$-category~$C(p)$
such that $C(p)\subset C(p+1)$.
\end{proposition}

We will use a construction of $\nu I^n$ based on chain complexes
(\cite{Ste}, Example~3.10). Recall that if $K$~and~$L$ are
augmented chain complexes, then their tensor product augmented
chain complex $K\otimes L$ is given by
\begin{align*}
 &(K\otimes L)_q=\bigoplus_{i+j=q}(K_i\otimes L_j),\\
 &\partial(x\otimes y)=\partial x\otimes y+(-1)^i x\otimes\partial y
 \ \text{if $x\in K_i$},\\
 &\epsilon(x\otimes y)=(\epsilon x)(\epsilon y).
\end{align*}
It is convenient to identify cubes with their chain complexes, as
follows.

\begin{definition} \label{3B}
The \emph{standard interval} is the free augmented chain
complex~$I$ concentrated in degrees $0$~and~$1$ such that
$I_1$~has basis~$u_1$, such that $I_0$~has basis $\partial^-
u_1,\partial^+ u_1$, such that $\partial u_1=\partial^+
u_1-\partial^- u_1$, and such that $\epsilon\partial^-
u_1=\epsilon\partial^+ u_1=1$. The \emph{standard $n$-cube}~$I^n$
is the $n$-fold tensor product $I\otimes\ldots\otimes I$. The
$n$-chain $u_1\otimes\ldots\otimes u_1$ in~$I^n$ is denoted~$u_n$.
\end{definition}

Thus the chain complex~$I^n$ is a free chain complex on a basis
got by taking $n$-fold tensor products of the basis elements
$u_1,\partial^- u_1,\partial^+ u_1$ for~$I$; in particular
$I^0$~is free on a single zero-dimensional basis element~$u_0$.
The elements of the tensor product basis for~$I^n$ are called
standard basis elements. They correspond to the faces of a
geometrical $n$-cube, or, equivalently, they correspond to the
precubical operations. In fact there is an obvious embedding of
the precubical category in the category of chain complexes as
follows.

\begin{proposition} \label{3Z}
There is a functor from the precubical category to the category of
chain complexes given on objects by $[n]\mapsto I^n$ and on
morphisms by
$$\check\partial_i^\alpha(x\otimes y)
 =x\otimes\partial^\alpha u_1\otimes y\
 \text{for $x\in I^{i-1}$ and $y\in I^{n-i}$}.$$
The morphisms in the image of this functor are
augmentation-preserving and take standard basis elements to
standard basis elements.
\end{proposition}

We will also need subcomplexes corresponding to shells and boxes,
as follows.

\begin{definition} \label{3C}
The \emph{standard $n$-shell}~$S^n$ is the subcomplex of~$I^n$
generated by the standard basis elements other than~$u_n$. If
$\sigma$~is an $(n-1)$-dimensional basis element for~$I^n$, then
the \emph{standard $n$-box}~$B(\sigma)$ opposite~$\sigma$ is the
subcomplex of~$I^n$ generated by the standard basis elements other
than $u_n$~and~$\sigma$.
\end{definition}

We give partial orderings to the chain groups in the standard
chain complexes by the rule that $x\geq y$ if and only if $x-y$ is
a sum of standard basis elements. We then associate
$\omega$-categories to these complexes as follows.

\begin{definition} \label{3D}
Let $K$ be a standard cube, shell or box. The \emph{associated
$\omega$-category} $\nu K$ is the set of double sequences
$$x=(x_0^-,x_0^+\mid x_1^-,x_1^+\mid\ldots\,),$$
where $x_i^-$~and~$x_i^+$ are $i$-chains in~$K$ such that
\begin{align*}
 &x_i^-\geq 0,\\
 &x_i^+\geq 0,\\
 &\epsilon x_0^-=\epsilon x_0^+=1,\\
 &x_i^+ -x_i^-=\partial x_{i+1}^-=\partial x_{i+1}^+.
\end{align*}
The left identity $d_p^- x$ and right identity $d_p^+ x$ of the
double sequence
$$x=(x_0^-,x_0^+\mid x_1^-,x_1^+\mid\ldots\,)$$
are given by
$$d_p^\alpha x=(x_0^-,x_0^+\mid \ldots
 \mid x_{p-1}^-,x_{p-1}^+\mid x_p^\alpha,x_p^\alpha\mid 0,0\mid
 \ldots\,).$$
If $d_p^+x=d_p^- y=w$, say, then the composite $x\comp_p y$ is
given by
$$x\comp_p y=x-w+y,$$
where the addition and subtraction are performed termwise.
\end{definition}

For example, in $\nu I^2$ there are elements
$$a=
 (\check\partial_2^-\check\partial_1^- u_0,
 \check\partial_2^+\check\partial_1^- u_0\mid
 \check\partial_1^- u_1,\check\partial_1^- u_1\mid
 0,0\mid\ldots\,)$$
and
$$b=
 (\check\partial_2^+\check\partial_1^- u_0,
 \check\partial_2^+\check\partial_1^+ u_0\mid
 \check\partial_2^+ u_1,\check\partial_2^+ u_1\mid
 0,0\mid\ldots\,)$$
such that
$$d_0^+ a
 =d_0^- b
 =(\check\partial_2^+\check\partial_1^- u_0,
 \check\partial_2^+\check\partial_1^- u_0\mid
 0,0\mid\ldots\,)$$
and
$$a\comp_0 b
 =(\check\partial_2^-\check\partial_1^- u_0,
 \check\partial_2^+\check\partial_1^+ u_0\mid
 \check\partial_1^- u_1+\check\partial_2^+ u_1,
 \check\partial_1^- u_1+\check\partial_2^+ u_1\mid
 0,0\mid\ldots\,).$$
The elements $a$, $b$ and $a\comp_0 b$ correspond to the morphisms
$\partial_1^- x$, $\partial_2^+ x$ and $\partial_1^-
x\circ\partial_2^+ x$ in Figure~\ref{F1}.

Given a chain~$x$ of positive degree we will write $\partial^- x$
and $\partial^+ x$ for the negative and positive parts of
$\partial x$; thus $\partial^- x$ and $\partial^+ x$ are the
linear combinations of disjoint families of basis elements with
positive integer coefficients such that $\partial x=\partial^+
x-\partial^- x$. This notation is of course consistent with the
earlier use of $\partial^- u_1$ and $\partial^+ u_1$. We can now
construct some specific double sequences of chains in~$I^n$, which
turn out to be members of $\nu I^n$, as follows.

\begin{definition} \label{3E}
Let $\sigma$ be a $p$-dimensional basis element for~$I^n$. Then
the associated \emph{atom}~$\langle\sigma\rangle$ is the double
sequence given by
\begin{align*}
 &\langle\sigma\rangle_q^\alpha=(\partial^\alpha)^{p-q}\sigma
 \ \text{for $q\leq p$},\\
 &\langle\sigma\rangle_q^\alpha=0
 \ \text{for $q>p$}.
\end{align*}
\end{definition}

\begin{figure}
 $$\xymatrix@!C{
 &&&\ar@{-}[]+0;[rrrr]+0^{\langle u_3\rangle_1^-}
 &&&&\ar[]+0;[ddrr]
 \\
 \\
 &\ar@{}[l]|<<<<{\langle u_3\rangle_0^-}
 \ar@{-}[uurr]+0
 \ar@{}[uurrrrrr]|{\check\partial_1^- u_2}
 \ar@{}[ddrrrrrr]|{\check\partial_3^- u_2}
 \ar@{-}[ddrr]+0
 &&&&\ar@{.}[uurr]
 \ar@{.}[llll]|<<<<{\langle u_3\rangle_2^-\Downarrow}
 \ar@{}[rrrr]|{\check\partial_2^+ u_2}
 \ar@{.}[ddrr]
 &&&&\ar@{}[r]|<<<<{\langle u_3\rangle_0^+}
 &\\
 \\
 &&&\ar@{-}[]+0;[rrrr]+0_{\langle u_3\rangle_1^+}
 &&&&\ar[]+0;[uurr]
 \\
 &&&&&\ar@3{->}[d]_{\langle u_3\rangle_3^-=\langle u_3\rangle_3^+}
 \\
 &&&&&\\
 &&&\ar@{-}[]+0;[rrrr]+0^{\langle u_3\rangle_1^-}
 \ar@{}[ddrrrrrr]|{\check\partial_3^+ u_2}
 &&&&\ar[]+0;[ddrr]
 \\
 \\
 &\ar@{}[l]|<<<<{\langle u_3\rangle_0^-}
 \ar@{-}[uurr]+0
 \ar@{}[rrrr]|{\check\partial_2^- u_2}
 \ar@{-}[ddrr]+0
 &&&&\ar@{.}[uull]
 \ar@{.}[rrrr]|<<<<{\langle u_3\rangle_2^+\Downarrow}
 \ar@{.}[ddll]
 &&&&\ar@{}[r]|<<<<{\langle u_3\rangle_0^+}
 &\\
 \\
 &&&\ar@{}[uurrrrrr]|{\check\partial_1^+ u_2}
 \ar@{-}[]+0;[rrrr]+0_{\langle u_3\rangle_1^+}
 &&&&\ar[]+0;[uurr]
 }$$
 \caption{The chains $\langle u_3\rangle_i^\alpha$}
 \label{F4}
\end{figure}
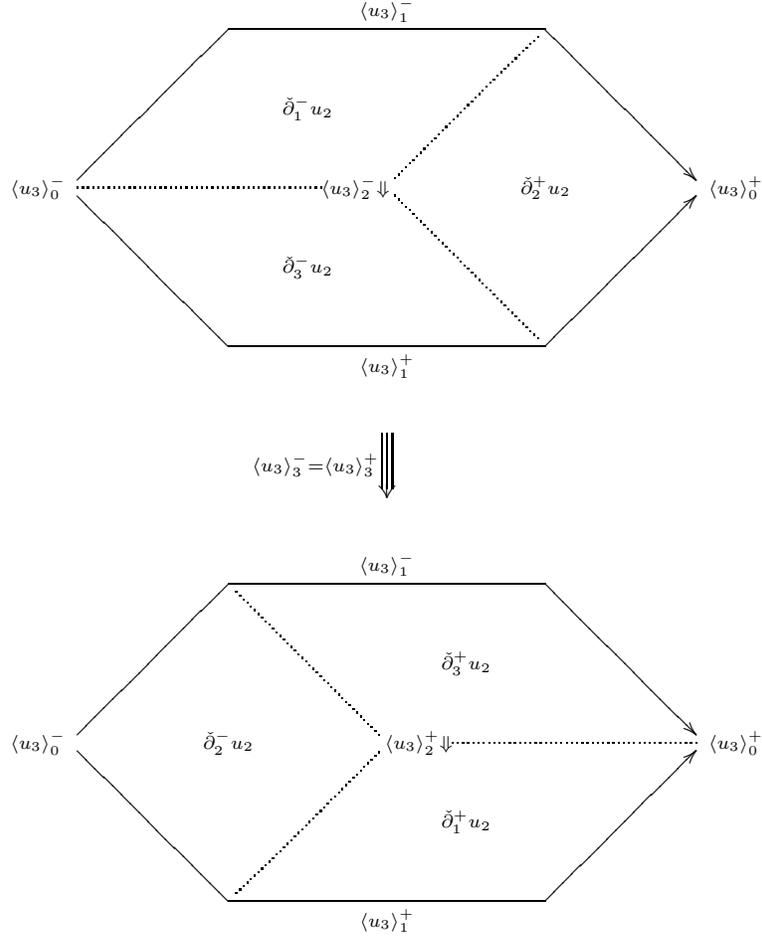

For example
\begin{align*}
 &\langle u_3\rangle_0^-
 =\check\partial_3^-\check\partial_2^-\check\partial_1^- u_0,\\
 &\langle u_3\rangle_0^+
 =\check\partial_3^+\check\partial_2^+\check\partial_1^+ u_0,\\
 &\langle u_3\rangle_1^-
 =\check\partial_2^-\check\partial_1^- u_1
 +\check\partial_3^+\check\partial_1^- u_1
 +\check\partial_3^+\check\partial_2^+ u_1,\\
 &\langle u_3\rangle_1^+
 =\check\partial_3^-\check\partial_2^- u_1
 +\check\partial_3^-\check\partial_1^+ u_1
 +\check\partial_2^+\check\partial_1^+ u_1,\\
 &\langle u_3\rangle_2^-
 =\check\partial_1^- u_2
 +\check\partial_2^+ u_2
 +\check\partial_3^- u_2,\\
 &\langle u_3\rangle_2^+
 =\check\partial_3^+ u_2
 +\check\partial_2^- u_2
 +\check\partial_1^+ u_2,\\
 &\langle u_3\rangle_3^-
 =\langle u_3\rangle_3^+
 =u_3;
\end{align*}
thus the terms of the chains~$\langle u_3\rangle_i^\alpha$
correspond to the precubical operations on $3$-cubes with standard
decompositions
$\partial_{i(1)}^{\alpha(1)}\ldots\partial_{i(p)}^{\alpha(p)}$
such that the signs $(-)^{i(r)-r}\alpha(r)$ are constant. These
are the precubical operations which were called extreme in
Section~\ref{S2}, and one can draw a $3$-cube with the chains
of~$\langle u_3\rangle$ on its extremities; see Figure~\ref{F4}.

Figure~\ref{F4} also illustrates complementary operations. For
example, the operations complementary to~$\partial_1^-$ are $\id$,
$\partial_2^+$, $\partial_3^-$, $\partial_2^+\partial_3^+$. In
Figure~\ref{F4} the face corresponding to $\id$ (the $3$-cube
itself) joins $\check\partial_1^- u_2$ to~$\langle
u_3\rangle_2^+$; the faces corresponding to
$\partial_2^+$~and~$\partial_3^-$ join $\check\partial_1^- u_2$
to~$\langle u_3\rangle_1^+$; the face corresponding to
$\partial_2^+\partial_3^+$ (the terminal edge in~$\langle
u_3\rangle_1^-$) joins $\check\partial_1^- u_2$ to $\langle
u_3\rangle_0^+$.

The formulae for the chains~$\langle u_3\rangle_i^\alpha$ extend
to all dimensions. For atoms in general, using the definition of
the boundary in a tensor product of chain complexes, we get the
following formulae.

\begin{proposition} \label{3F}
The atoms in $\nu I$ are given by
\begin{align*}
 &\langle u_1\rangle_0^\alpha=\partial^\alpha u_1,\\
 &\langle u_1\rangle_1^\alpha=u_1,\\
 &\langle u_1\rangle_q^\alpha=0\ \text{for $q>1$},\\
 &\langle\partial^\beta u_1\rangle_0^\alpha=\partial^\beta u_1,\\
 &\langle\partial^\beta u_1\rangle_q^\alpha=0\ \text{for $q>0$}.
\end{align*}
The atoms in $\nu I^n$ are given by
$$\langle\sigma_1\otimes\ldots\otimes\sigma_n\rangle_q^\alpha
 =\sum_{i(1)+\ldots+i(n)=q}
 \langle\sigma_1\rangle_{i(1)}^\alpha
 \otimes\langle\sigma_2\rangle_{i(2)}^{(-)^{i(1)}\alpha}
 \otimes\ldots
 \otimes\langle\sigma_n\rangle_{i(n)}^{(-)^{i(1)+\ldots+i(n-1)}\alpha}.$$
\end{proposition}

From Proposition~\ref{3F}, if $\sigma$~is a basis element in~$I^n$
then $\epsilon\langle\sigma\rangle_0^\alpha=1$, and it follows
that the atoms belong to $\nu I^n$. In fact they generate $\nu
I^n$, and the standard $\omega$-categories have presentations in
terms of these generators as follows (\cite{Ste}, Theorem~6.1).

\begin{theorem} \label{3H}
Let $K$ be a standard cube, shell or box. Then $\nu K$ is freely
generated by the atoms corresponding to the basis elements of~$K$
subject to the following relations\textup{:} if
$\langle\sigma\rangle$ is an atom corresponding to a
$p$-dimensional basis element~$\sigma$, then
$$d_p^-\langle\sigma\rangle
 =d_p^+\langle\sigma\rangle
 =\langle\sigma\rangle;$$
if $\langle\sigma\rangle$ is an atom corresponding to a
$p$-dimensional basis element~$\sigma$ with $p>0$, then
$d_{p-1}^-\langle\sigma\rangle=w^-$ and
$d_{p-1}^+\langle\sigma\rangle=w^+$, where $w^-$~and~$w^+$ are
suitable composites of lower-dimensional atoms.
\end{theorem}

For example, the presentation of $\nu I^2$ is essentially as
illustrated in Figure~\ref{F1}. There are generators
$\langle\check\partial_2^\beta\check\partial_1^\alpha u_0\rangle$
subject to relations
$d_0^\gamma\langle\check\partial_2^\beta\check\partial_1^\alpha
u_0\rangle=\langle\check\partial_2^\beta\check\partial_1^\alpha
u_0\rangle$; there are generators $\langle\check\partial_1^\alpha
u_1\rangle$ and $\langle\check\partial_2^\alpha u_1\rangle$
subject to relations
$$d_1^\gamma\langle\check\partial_i^\alpha
 u_1\rangle=\langle\check\partial_i^\alpha u_1\rangle,\
 d_0^\gamma\langle\check\partial_1^\alpha
 u_1\rangle=\langle\check\partial_2^\gamma\check\partial_1^\alpha
 u_0\rangle,\ d_0^\gamma\langle\check\partial_2^\alpha
 u_1\rangle=\langle\check\partial_2^\alpha\check\partial_1^\gamma
 u_0\rangle;$$ there is a generator~$\langle u_2\rangle$ subject to
relations
$$d_2^\gamma\langle u_2\rangle=\langle u_2\rangle,\
 d_1^-\langle u_2\rangle
  =\langle\check\partial_1^-u_1\rangle
  \comp_0\langle\check\partial_2^+u_1\rangle,\
 d_1^+\langle u_2\rangle
  =\langle\check\partial_2^-u_1\rangle
  \comp_0\langle\check\partial_1^+u_1\rangle.$$

A more detailed description of the way in which the atoms generate
$\nu I^n$ is as follows (\cite{Ste}, Proposition~5.4).

\begin{theorem} \label{3I}
Let
$$x=(x_0^-,x_0^+\mid x_1^-,x_1^+\mid\ldots\,)$$
be an element of $\nu I^n$. Then $x$~is an identity for~$\comp_p$
if and only if $x_q^-=x_q^+=0$ for $q>p$. If $x$~is an identity
for~$\comp_p$, then
$$x_p^-=x_p^+=\sigma_1+\ldots+\sigma_k$$
for some $p$-dimensional basis elements~$\sigma_i$, and $x$~is a
composite of the atoms~$\langle\sigma_i\rangle$ together with
atoms of lower dimension.
\end{theorem}

\begin{figure}
 $$\xymatrix@!C{
 &&\ar@{-}[]+0;[rrrr]+0&&&&
 \ar@{-}[]+0;[ddrr]+0
 ^{\langle\check\partial_3^+\check\partial_2^+u_1\rangle}\\
 &&&&&&\ar@{-}[]+0;[ddrr]+0
 \\
 \ar@{-}[]+0;[uurr]+0
 \ar@{}[uurrrrrr]|{\langle\check\partial_1^- u_2\rangle}
 \ar@{-}[]+0;[rrrr]+0&&&&
 \ar@{-}[]+0;[uurr]+0&&&&\\
 \ar@{-}[]+0;[rrrr]+0
 ^{\langle\check\partial_3^-\check\partial_1^-u_1\rangle}&&&&
 \ar@{-}[]+0;[uurr]+0
 \ar@{}[rrrr]|{\langle\check\partial_2^+ u_2\rangle}
 \ar@{-}[]+0;[ddrr]+0&&&&\\
 \ar@{-}[]+0;[rrrr]+0
 \ar@{}[ddrrrrrr]|{\langle\check\partial_3^- u_2\rangle}
 \ar@{-}[]+0;[ddrr]+0&&&&
 \ar@{-}[]+0;[ddrr]+0&&&&\\
 &&&&&&\ar@{-}[]+0;[uurr]+0\\
 &&\ar@{-}[]+0;[rrrr]+0&&&&
 \ar@{-}[]+0;[uurr]+0
 _{\langle\check\partial_2^+\check\partial_1^+u_1\rangle}
 }$$
 \caption{The decomposition of $d_2^-\langle u_3\rangle$}
 \label{F5}
\end{figure}
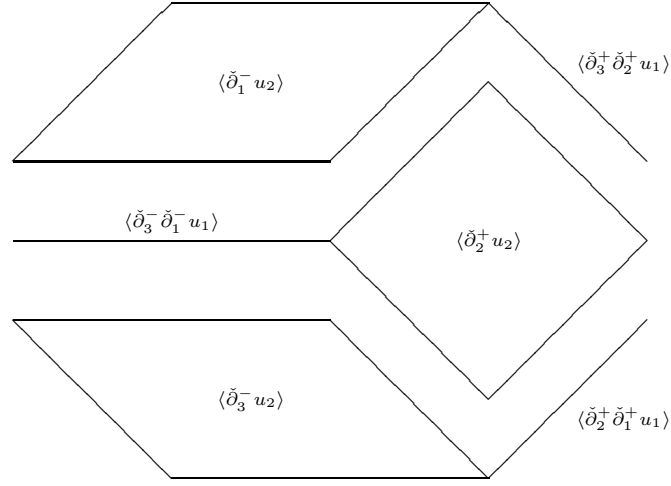

\begin{example} \label{3Y}
In $\nu I^3$ the element
$$d_2^-\langle u_3\rangle
 =(\langle u_3\rangle_0^-,\langle u_3\rangle_0^+\mid
 \langle u_3\rangle_1^-,\langle u_3\rangle_1^+\mid
 \langle u_3\rangle_2^-,\langle u_3\rangle_2^-\mid
 0,0\mid\ldots\,),$$
is an identity for~$\comp_2$ with $2$-chain component
$$\langle u_3\rangle_2^-
 =\check\partial_1^- u_2
 +\check\partial_2^+ u_2
 +\check\partial_3^- u_2$$
and with decomposition
$$(\langle\check\partial_1^- u_2\rangle
 \comp_0\langle\check\partial_3^+\check\partial_2^+ u_1\rangle)
 \comp_1
 (\langle\check\partial_3^-\check\partial_1^- u_1\rangle
 \comp_0\langle\check\partial_2^+ u_2\rangle)
 \comp_1
 (\langle\check\partial_3^- u_2\rangle
 \comp_0\langle\check\partial_2^+\check\partial_1^+ u_1\rangle).$$
There is a corresponding decomposition of the top half of
Figure~\ref{F4}, which is shown in Figure~\ref{F5}. Similarly,
there is a decomposition of $d_2^+\langle u_3\rangle$ given by
$$(\langle\check\partial_2^-\check\partial_1^- u_1\rangle
 \comp_0\langle\check\partial_3^+ u_2\rangle)
 \comp_1
 (\langle\check\partial_2^- u_2\rangle
 \comp_0\langle\check\partial_3^+\check\partial_1^+ u_1\rangle)
 \comp_1
 (\langle\check\partial_3^-\check\partial_2^- u_1\rangle
 \comp_0\langle\check\partial_1^+ u_2\rangle).$$
\end{example}

Recall the functor from the precubical category to the category of
chain complexes given in Proposition~\ref{3Z}. Since the morphisms
in the image are augmentation-preserving and take standard basis
elements to standard basis elements, they induce morphisms between
the $\omega$-categories $\nu I^n$; thus there is a functor
$[n]\mapsto\nu I^n$ from the precubical category to the category
of $\omega$-categories. The functors $\hom(\nu I^n,-)$ therefore
produce a functor from $\omega$-categories to precubical sets. In
fact we get a functor from $\omega$-categories to stratified
precubical sets, as follows.

\begin{definition} \label{3J}
The \emph{cubical nerve} of an $\omega$-category~$C$ is the
stratified precubical set~$X$ given by
$$X_n=\hom(\nu I^n,C),$$
where an $n$-cube $x\colon\nu I^n\to C$ with $n>0$ is thin if
$x\langle u_n\rangle$ is an identity for~$\comp_{n-1}$.
\end{definition}

From the structure of the $\omega$-categories associated to the
standard complexes (Theorem~\ref{3H}), one sees that $\nu S^n$ is
generated by $2n$~copies of $\nu I^{n-1}$ corresponding to the
$(n-1)$-dimensional faces of an $n$-cube, subject to relations
corresponding to their pairwise intersections, and similarly for
boxes. This gives the following result.

\begin{proposition} \label{3K}
Let $C$ be an $\omega$-category. Then the $n$-shells in the
cubical nerve of~$C$ correspond to the morphisms $\nu S^n\to C$,
and the $n$-boxes opposite~$\partial_k^\gamma$ correspond to the
morphisms $\nu B(\check\partial_k^\gamma u_{n-1})\to C$. A filler
for an $n$-shell~$s$ is an $n$-cube~$x$ such that $x|\nu S^n=s$,
and similarly for boxes.
\end{proposition}

It remains to show that the cubical nerve of an $\omega$-category
has thin fillers satisfying the conditions for a complete
stratification. For $n>0$ it follows from Theorem~\ref{3H} that
$\nu I^n$ is got from $\nu S^n$ by adjoining a single extra
generator~$\langle u_n\rangle$ subject to relations on the
$d_n^\alpha\langle u_n\rangle$ and the $d_{n-1}^\alpha\langle
u_n\rangle$. This gives us the following result on thin fillers
for shells.

\begin{proposition} \label{3L}
Let $s$ be an $n$-shell with $n>0$ in the cubical nerve of an
$\omega$-category. Then $s$~has a thin filler if and only if
$s(d_{n-1}^-\langle u_n\rangle)=s(d_{n-1}^+\langle u_n\rangle)$.
If $s$~does have a thin filler, then this thin filler is unique.
\end{proposition}

\begin{proof}
From Theorem~\ref{3H}, the fillers~$x$ of~$s$ correspond to
elements $x\langle u_n\rangle$ such that
\begin{align*}
 &d_n^- x\langle u_n\rangle=d_n^+ x\langle u_n\rangle
 =x\langle u_n\rangle,\\
 &d_{n-1}^- x\langle u_n\rangle=s(d_{n-1}^-\langle u_n\rangle),\\
 &d_{n-1}^+ x\langle u_n\rangle=s(d_{n-1}^+\langle u_n\rangle).
\end{align*}
The filler~$x$ is thin if and only if $x\langle u_n\rangle$~is an
identity for~$\comp_{n-1}$; thus $x$~is thin if and only if it
satisfies the additional relations
$$d_{n-1}^- x\langle u_n\rangle=d_{n-1}^+ x\langle u_n\rangle
 =x\langle u_n\rangle.$$
These additional relations actually imply that $d_n^- x\langle
u_n\rangle=d_n^+ x\langle u_n\rangle=x\langle u_n\rangle$, because
$d_n^\alpha d_{n-1}^\alpha=d_{n-1}^\alpha$ (see
Proposition~\ref{3A}), so the defining relations for thin fillers
reduce to
$$x\langle u_n\rangle
 =s(d_{n-1}^-\langle u_n\rangle)
 =s(d_{n-1}^+\langle u_n\rangle).$$
It follows that $s$~has a unique thin filler if $s(d_{n-1}^-
u_n)=s(d_{n-1}^+ u_n)$, and that $s$~has no thin filler otherwise.
This completes the proof.
\end{proof}

In Section~\ref{S6} we will prove the following result.

\begin{theorem} \label{3O}
Let $\partial_k^\gamma$ be a face operation on $n$-cubes. Then
there is a factorisation
\begin{align*}
 &d_{n-1}^{(-)^{k-1}\gamma}\langle u_n\rangle\\
 &=A_{n-1}^-\comp_{n-2}(A_{n-2}^-\comp_{n-3}\ldots (A_1^-\comp_0
  \langle\check\partial_k^\gamma u_{n-1}\rangle
  \comp_0 A_1^+)\ldots\comp_{n-3}A_{n-2}^+)\comp_{n-2}A_{n-1}^+
\end{align*}
in $\nu I^n$ such that the~$A_q^\alpha$ are in $\nu
B(\check\partial_k^\gamma u_{n-1})$ and such that $b(A_q^-)$ and
$b(A_q^+)$ are identities for~$\comp_{q-1}$ whenever $b$~is an
admissible $n$-box opposite~$\partial_k^\gamma$.
\end{theorem}

For example, consider the operation~$\partial_2^+$ on $3$-cubes.
From Example~\ref{3Y} we get
$$d_2^-\langle u_3\rangle
 =A_2^-\comp_1
 (A_1^-\comp_0\langle\check\partial_2^+ u_2\rangle\comp_0 A_1^+)
 \comp_1 A_2^+,$$
where
\begin{align*}
 &A_2^-=\langle\check\partial_1^- u_2\rangle
 \comp_0\langle\check\partial_3^+\check\partial_2^+ u_1\rangle,\\
 &A_2^+=\langle\check\partial_3^- u_2\rangle
 \comp_0\langle\check\partial_2^+\check\partial_1^+ u_1\rangle,\\
 &A_1^-=\langle\check\partial_3^-\check\partial_1^- u_1\rangle,
\end{align*}
and $A_1^+$~is an identity for~$\comp_0$. We see that the atomic
factors of~$A_q^\alpha$ have dimension at most~$q$; we also see
that the $q$-dimensional factors of the~$A_q^\alpha$ correspond to
the precubical operations $\partial_1^-$, $\partial_3^-$ and
$\partial_1^-\partial_3^-$, which are the non-identity precubical
operations complementary to~$\partial_2^+$. If $b$~is an
admissible $3$-box opposite~$\partial_2^+$, then $\partial_1^- b$,
$\partial_3^- b$ and $\partial_1^-\partial_3^- b$ are thin, so
$b\langle\tau\rangle$ is an identity for~$\comp_{q-1}$ whenever
$\langle\tau\rangle$~is an atomic factor in~$A_q^\alpha$, and it
follows from Proposition~\ref{3A} that $b(A_q^\alpha)$ is an
identity for~$\comp_{q-1}$ as claimed.

Assuming Theorem~\ref{3O} we get the following result on
admissible boxes.

\begin{theorem} \label{3M}
Let $b$ be an admissible $n$-box opposite~$\partial_k^\gamma$. If
$s$~is an $n$-shell extending~$b$ then
$$s\langle\check\partial_k^\gamma u_{n-1}\rangle
 =s(d_{n-1}^{(-)^{k-1}\gamma}\langle u_n\rangle).$$
If $n\geq 2$ then
$$b(d_{n-2}^\alpha\langle\check\partial_k^\gamma u_{n-1}\rangle)
 =d_{n-2}^\alpha b(d_{n-1}^{(-)^k\gamma}\langle u_n\rangle).$$
\end{theorem}

(Note here that $b(d_{n-1}^{(-)^k\gamma}\langle u_n\rangle)$
exists because, by Proposition~\ref{3F}, $\check\partial_k^\gamma
u_{n-1}$ is not a term in~$\langle
u_n\rangle_{n-1}^{(-)^k\gamma}$.)

\begin{proof}
Let $s$ be an $n$-shell extending~$b$ and apply~$s$ to the
factorisation in Theorem~\ref{3O}. Since
$s(A_q^\beta)=b(A_q^\beta)$ is an identity for~$\comp_{q-1}$, the
factorisation collapses to the equality
$s\langle\check\partial_k^\gamma
u_{n-1}\rangle=s(d_{n-1}^{(-)^{k-1}\gamma}\langle u_n\rangle)$.

Now suppose that $n\geq 2$. By Proposition~\ref{3A} and a collapse
like that in the previous paragraph,
$$b(d_{n-2}^\alpha\langle\check\partial_k^\gamma u_{n-1}\rangle)
 =b(d_{n-2}^\alpha d_{n-1}^{(-)^{k-1}\gamma}\langle u_n\rangle).$$
By a further application of Proposition~\ref{3A},
$$b(d_{n-2}^\alpha d_{n-1}^{(-)^{k-1}\gamma}\langle u_n\rangle)
 =b(d_{n-2}^\alpha d_{n-1}^{(-)^k\gamma}\langle u_n\rangle)
 =d_{n-2}^\alpha b(d_{n-1}^{(-)^k\gamma}\langle u_n\rangle);$$
therefore $b(d_{n-2}^\alpha\langle\check\partial_k^\gamma
u_{n-1}\rangle)=d_{n-2}^\alpha b(d_{n-1}^{(-)^k\gamma}\langle
u_n\rangle)$ as required. This completes the proof.
\end{proof}

We now get the main theorem of this section as follows.

\begin{theorem} \label{3N}
The cubical nerve of a strict $\omega$-category is a complete
stratified precubical set.
\end{theorem}

\begin{proof}
We will verify the conditions of Definition~\ref{2M}.

Let $s$ be an admissible $n$-shell; we must show that $s$~has a
unique thin filler. By the definition of an admissible shell there
are distinct non-complementary face operations
$\partial_k^\gamma$~and~$\partial_l^\delta$ such that
$s_k^\gamma=s_l^\delta$ and such that the boxes formed by removing
$s_k^\gamma$~or~$s_l^\delta$ are admissible. Since
$\partial_k^\gamma$~and~$\partial_l^\delta$ are not complementary
we have $(-)^k\gamma\neq (-)^l\delta$ (see Remark~\ref{2G}). Of
the expressions $s(d_{n-1}^-\langle u_n\rangle)$ and
$s(d_{n-1}^+\langle u_n\rangle)$ it follows from Theorem~\ref{3M}
that one is equal to $s\langle\check\partial_k^\gamma
u_{n-1}\rangle$ and the other is equal to
$s\langle\check\partial_l^\delta u_{n-1}\rangle$; in other words,
one of them is $s_k^\gamma\langle u_{n-1}\rangle$ and the other is
$s_l^\delta\langle u_{n-1}\rangle$. But $s_k^\gamma=s_l^\delta$,
so $s(d_{n-1}^-\langle u_n\rangle)=s(d_{n-1}^+\langle
u_n\rangle)$. By Proposition~\ref{3L}, $s$~has a unique thin
filler.

Now let $b$ be an admissible $n$-box opposite~$\partial_k^\gamma$;
we must show that $b$~has a unique thin filler. Because of
Proposition~\ref{3L}, the thin fillers of~$b$ correspond to
$n$-shells~$s$ extending~$b$ such that
$$s(d_{n-1}^-\langle u_n\rangle)=s(d_{n-1}^+\langle u_n\rangle).$$
Equivalently, by Theorem~\ref{3M}, these are the $n$-shells~$s$
extending~$b$ such that
$$s\langle\check\partial_k^\gamma u_{n-1}\rangle
 =b(d_{n-1}^{(-)^k\gamma}\langle u_n\rangle).$$
Now it follows from Theorem~\ref{3H} that $\nu S^n$ is got from
$\nu B(\check\partial_k^\gamma u_{n-1})$ by adjoining
$\langle\check\partial_k^\gamma u_{n-1}\rangle$ and imposing
certain relations. An $n$-shell~$s$ extending~$b$ is therefore
uniquely determined by the value of
$s\langle\check\partial_k^\gamma u_{n-1}\rangle$, and the possible
values for $s\langle\check\partial_k^\gamma u_{n-1}\rangle$ are
given by the following conditions: in all cases, we require
$$d_{n-1}^-s\langle\check\partial_k^\gamma u_{n-1}\rangle
 =d_{n-1}^+s\langle\check\partial_k^\gamma u_{n-1}\rangle
 =s\langle\check\partial_k^\gamma u_{n-1}\rangle;$$
if $n\geq 2$ then we also require
$$d_{n-2}^\alpha s\langle\check\partial_k^\gamma u_{n-1}\rangle
 =b(d_{n-2}^\alpha\langle\check\partial_k^\gamma u_{n-1}\rangle).$$
Using Theorem~\ref{3M} in the case $n\geq 2$, we see that these
conditions are satisfied when $s\langle\check\partial_k^\gamma
u_{n-1}\rangle=b(d_{n-1}^{(-)^k\gamma}\langle u_n\rangle)$.
Therefore $b$~has a unique thin filler.

Finally let $b$ be an admissible $n$-box
opposite~$\partial_k^\gamma$ such that all the
$(n-1)$-cubes~$b_i^\alpha$ are thin, and let $x$ be the thin
filler of~$b$; we must show that $\partial_k^\gamma x$ is thin.
Since $x$~is thin we have
$$x(d_{n-1}^-\langle u_n\rangle)
 =d_{n-1}^- x\langle u_n\rangle
 =x\langle u_n\rangle
 =d_{n-1}^+ x\langle u_n\rangle
 =x(d_{n-1}^+ \langle u_n\rangle),$$
and it then follows from Theorem~\ref{3M} that
$$(\partial_k^\gamma x)\langle u_{n-1}\rangle
 =x(\check\partial_k^\gamma\langle u_{n-1}\rangle)
 =x(d_{n-1}^{(-)^{k-1}\gamma}\langle u_n\rangle)
 =x(d_{n-1}^{(-)^k\gamma}\langle u_n\rangle)
 =b(d_{n-1}^{(-)^k\gamma}\langle u_n\rangle).$$
Since the~$b_i^\alpha$ are thin, $b\langle\tau\rangle$ is an
identity for~$\comp_{n-2}$ whenever $\langle\tau\rangle$~is an
atom in $\nu B(\check\partial_k^\gamma u_{n-1})$. By
Theorem~\ref{3H}, these atoms generate $\nu
B(\check\partial_k^\gamma u_{n-1})$, so
$b(d_{n-1}^{(-)^k\gamma}\langle u_n\rangle)$ is an identity
for~$\comp_{n-2}$ by Proposition~\ref{3A}. Therefore
$(\partial_k^\gamma x)\langle u_{n-1}\rangle$ is an identity
for~$\comp_{n-2}$, which means that $\partial_k^\gamma x$ is thin.

This completes the proof.
\end{proof}

\section{From complete stratified precubical sets to omega-categories}
\label{S4}

Throughout this section, let $X$ be a complete stratified
precubical set. We will show that $X$~is the cubical nerve of an
$\omega$-category by constructing degeneracies, connections and
compositions with the properties of \cite{ABS}.

The degeneracies are to be operations
$$\epsilon_1,\ldots,\epsilon_n\colon X_{n-1}\to X_n,$$
and we will define $\epsilon_k x$ for $x\in X_{n-1}$ as the unique
thin filler of an admissible $n$-shell; in other words,
$\epsilon_k x$ is a thin $n$-cube with prescribed values for the
faces $\partial_i^\alpha\epsilon_k x$. The process is inductive:
we define degeneracies on $n$-cubes in terms of degeneracies on
lower-dimensional cubes.

\begin{definition} \label{4A}
The \emph{degeneracies} are the elements $\epsilon_k x$, defined
for $k=1$, $2$, \dots,~$n$ and $x\in X_{n-1}$, such that
$\epsilon_k x$ is a thin member of~$X_n$ and
\begin{align*}
 &\partial_i^\alpha\epsilon_k x
  =\epsilon_{k-1}\partial_i^\alpha x\ \text{for $i<k$},\\
 &\partial_k^\alpha\epsilon_k x=x,\\
 &\partial_i^\alpha\epsilon_k x=\epsilon_k\partial_{i-1}^\alpha x
  \ \text{for $i>k$}.
\end{align*}
\end{definition}

\begin{figure}
 $$\xymatrix@!C{
 \ar@{-}[]+0;[rr]+0^{\epsilon_1\partial_1^+ x}_{=}&&&&
 \ar@{-}[]+0;[rr]+0^{x}&&\\
 &\epsilon_1 x&&&&\epsilon_2 x\\
 \ar@{-}[]+0;[uu]+0^{x}
 \ar@{-}[]+0;[rr]+0_{\epsilon_1\partial_1^- x}^{=}&&
 \ar@{-}[]+0;[uu]+0_{x}&&
 \ar@{-}[]+0;[uu]+0^{\epsilon_1\partial_1^- x}_{=}
 \ar@{-}[]+0;[rr]+0_{x}&&
 \ar@{-}[]+0;[uu]+0_{\epsilon_1\partial_1^+ x}^{=}
 }$$
 \caption{The degeneracies of a $1$-cube}
 \label{F6}
\end{figure}
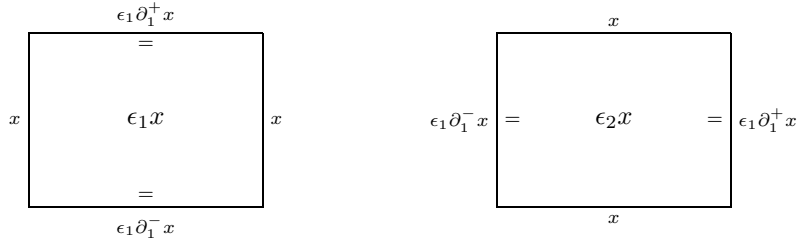

The two degeneracies of a $1$-cube~$x$ are shown in
Figure~\ref{F6}, with the thin edges labelled by equality signs.
Compare the first two shells in Figure~\ref{F3}.

To justify Definition~\ref{4A}, we must show that the prescribed
values~$s_i^\alpha$ for $\partial_i^\alpha\epsilon_k x$ form an
admissible $n$-shell~$s$. We will in fact show that $s$~is of the
type described in Example~\ref{2K}, using induction on~$n$.
Suppose that there are degeneracies with the required properties
on $m$-cubes for $m<n-1$. Then we get $\partial_i^\alpha
s_j^\beta=\partial_j^\beta s_{i+1}^\alpha$ for $i\geq j$ as
follows: if $j\leq i<k-1$ then
$$\partial_i^\alpha s_j^\beta
 =\partial_i^\alpha\epsilon_{k-1}\partial_j^\beta x
 =\epsilon_{k-2}\partial_i^\alpha\partial_j^\beta x
 =\epsilon_{k-2}\partial_j^\beta\partial_{i+1}^\alpha x
 =\partial_j^\beta\epsilon_{k-1}\partial_{i+1}^\alpha x
 =\partial_j^\beta s_{i+1}^\alpha,$$
if $j\leq i=k-1$ then
$$\partial_i^\alpha s_j^\beta
 =\partial_i^\alpha\epsilon_{k-1}\partial_j^\beta x
 =\partial_j^\beta x
 =\partial_j^\beta s_{i+1}^\alpha,$$
etc. Therefore $s$~is a shell. We also have $s_k^-=s_k^+$. If
$\theta$~is a non-identity precubical operation on $n$-cubes whose
standard decomposition does not contain
$\partial_k^-$~or~$\partial_k^+$, say
$$\theta
 =(\partial_{i(1)}^{\alpha(1)}\ldots\partial_{i(p)}^{\alpha(p)})
 (\partial_{j(1)}^{\beta(1)}\ldots\partial_{j(q)}^{\beta(q)})$$
with $p+q>0$ and with
$$i(1)<\ldots<i(p)<k<j(1)<\ldots<j(q),$$
then
$$\theta s
 =\epsilon_{k-p}
 (\partial_{i(1)}^{\alpha(1)}\ldots\partial_{i(p)}^{\alpha(p)})
 (\partial_{j(1)-1}^{\beta(1)}\ldots\partial_{j(q)-1}^{\beta(q)})x,$$
so $\theta s$ is thin. It follows from Example~\ref{2K} that
$s$~is an admissible shell, as required.

Connections are defined by a similar inductive process, with a
similar inductive justification using Example~\ref{2L}.

\begin{definition} \label{4B}
The \emph{connections} are the elements $\Gamma_k^- x$ and
$\Gamma_k^+ x$, defined for $k=1$, $2$, \dots,~$n$ and $x\in X_n$,
such that $\Gamma_k^\gamma x$ is a thin member of~$X_{n+1}$ and
\begin{align*}
 &\partial_i^\alpha\Gamma_k^\gamma x
  =\Gamma_{k-1}^\gamma\partial_i^\alpha x\ \text{for $i<k$},\\
 &\partial_k^\gamma\Gamma_k^\gamma x
 =\partial_{k+1}^\gamma\Gamma_k^\gamma x
 =x,\\
 &\partial_k^{-\gamma}\Gamma_k^\gamma x
 =\partial_{k+1}^{-\gamma}\Gamma_k^\gamma x
 =\epsilon_k\partial_k^{-\gamma}x,\\
 &\partial_i^\alpha\Gamma_k^\gamma x
 =\Gamma_k^\gamma\partial_{i-1}^\alpha x
  \ \text{for $i>k+1$}.
\end{align*}
\end{definition}

\begin{figure}
 $$\xymatrix{
 \ar@{-}[]+0;[rr]+0^{\epsilon_1\partial_1^+ x}_{=}&&&&&
 \ar@{-}[]+0;[rr]+0^{x}&&\\
 &\Gamma_1^- x&&&&&\Gamma_1^+ x\\
 \ar@{-}[]+0;[uu]+0^{x}
 \ar@{-}[]+0;[rr]+0_{x}&&
 \ar@{-}[]+0;[uu]+0_{\epsilon_1\partial_1^+ x}^{=}&&&
 \ar@{-}[]+0;[uu]+0^{\epsilon_1\partial_1^- x}_{=}
 \ar@{-}[]+0;[rr]+0_{\epsilon_1\partial_1^- x}^{=}&&
 \ar@{-}[]+0;[uu]+0_{x}
 }$$
 \caption{The connections of a $1$-cube}
 \label{F7}
\end{figure}
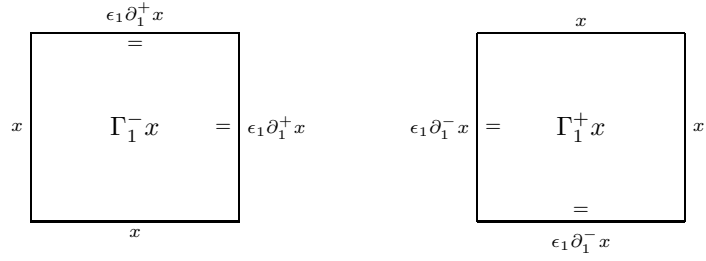

The two connections of a $1$-cube~$x$ are shown in
Figure~\ref{F7}; compare the last two shells in Figure~\ref{F3}.

We will get composites $x\circ_k y$ as the additional faces
$\partial_k^- G_k(x,y)$ of thin fillers of admissible boxes
opposite~$\partial_k^-$. These thin fillers are called composers,
and are again defined inductively. The definitions are as follows.

\begin{definition} \label{4C}
The \emph{composers} are the elements $G_k(x,y)$, defined for
$k=1$, $2$, \dots,~$n$ and $x,y\in X_n$ with
$\partial_k^+x=\partial_k^- y$, such that $G_k(x,y)$ is a thin
member of~$X_{n+1}$ and
\begin{align*}
 &\partial_i^\alpha G_k(x,y)
  =G_{k-1}(\partial_i^\alpha x,\partial_i^\alpha y)
  \ \text{for $i<k$},\\
 &\partial_k^+ G_k(x,y)=y,\\
 &\partial_{k+1}^- G_k(x,y)=x,\\
 &\partial_{k+1}^+ G_k(x,y)=\epsilon_k\partial_k^+ y,\\
 &\partial_i^\alpha G_k(x,y)
 =G_k(\partial_{i-1}^\alpha x,\partial_{i-1}^\alpha y)
  \ \text{for $i>k+1$}.
\end{align*}

The \emph{composites} are the elements
$$x\circ_k y=\partial_k^- G_k(x,y)\in X_n,$$
defined for $k=1$, $2$, \dots,~$n$ and $x,y\in X_n$ with
$\partial_k^+x=\partial_k^- y$.
\end{definition}

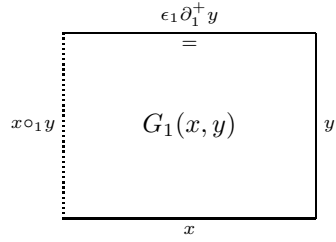
\begin{figure}
 $$\xymatrix{
 \ar@{-}[]+0;[rr]+0^{\epsilon_1\partial_1^+ y}_{=}&&\\
 &G_1(x,y)\\
 \ar@{.}[]+0;[uu]+0^{x\circ_1 y}
 \ar@{-}[]+0;[rr]+0_{x}&&
 \ar@{-}[]+0;[uu]+0_{y}
 }$$
 \caption{The composer of two $1$-cubes}
 \label{F8}
\end{figure}

The case $k=n=1$ is shown in Figure~\ref{F8}; compare
Figure~\ref{F2}. Note that $\partial_k^- G_k(x,y)$ is not
specified in the definition of $G_k(x,y)$, because we are dealing
with a box~$b$ opposite~$\partial_k^-$. The justification of this
definition is as before; in particular $\theta b$ is thin if
$\theta$~is a non-identity precubical operation whose standard
decomposition has no factors $\partial_k^-$, $\partial_k^+$
or~$\partial_{k+1}^-$, so $b$~is admissible by Example~\ref{2J}.

It remains to verify that $X$~satisfies the axioms for a cubical
nerve as given in \cite{ABS}. We begin with the following
observation.

\begin{proposition} \label{4E}
If $x\circ_k y$ is a composite such that $x$~and~$y$ are thin,
then $x\circ_k y$ is thin.
\end{proposition}

\begin{proof}
If $x$~and~$y$ are thin $n$-cubes then all the $n$-cubes in the
admissible box defining $G_k(x,y)$ are thin, so the additional
face $x\circ_k y=\partial_k^- G_k(x,y)$ is also thin.
\end{proof}

\begin{proposition} \label{4F}
The degeneracies have the property that
$$\epsilon_k\epsilon_l x=\epsilon_{l+1}\epsilon_k x
 \ \text {for $k\leq l$}.$$
\end{proposition}

\begin{proof}
The proof is by induction on~$n$, where $x$~is an $n$-cube. Using
the inductive hypothesis, we find that
$\partial_i^\alpha\epsilon_k\epsilon_l
x=\partial_i^\alpha\epsilon_{l+1}\epsilon_k x$ for
all~$\partial_i^\alpha$. This means that $\epsilon_k\epsilon_l x$
and $\epsilon_{l+1}\epsilon_k x$ are fillers for the same
shell~$s$, and in fact they are thin fillers for~$s$. But $s$~is
admissible (it is the admissible shell used to define
$\epsilon_k\epsilon_l x$), so it has a unique thin filler.
Therefore $\epsilon_k\epsilon_l x=\epsilon_{l+1}\epsilon_k x$.
\end{proof}

\begin{proposition} \label{4G}
The connections have the following properties\textup{:}
\begin{align*}
 &\Gamma_k^\gamma\epsilon_l x=\epsilon_{l+1}\Gamma_k^\gamma x
 \ \text{for $k<l$},\\
 &\Gamma_k^\gamma\epsilon_k x=\epsilon_{k+1}\epsilon_k x,\\
 &\Gamma_k^\gamma\epsilon_l x=\epsilon_l\Gamma_{k-1}^\gamma x
 \ \text{for $k>l$},\\
 &\Gamma_k^\gamma\Gamma_l^\delta x
 =\Gamma_{l+1}^\delta\Gamma_k^\gamma x
 \ \text{for $k<l$},\\
 &\Gamma_k^\gamma\Gamma_k^\gamma x
 =\Gamma_{k+1}^\gamma\Gamma_k^\gamma x.
\end{align*}
\end{proposition}

\begin{proof}
Similar.
\end{proof}

\begin{proposition} \label{4H}
If $x\circ_k y$ is defined, then
\begin{align*}
&\partial_i^\alpha(x\circ_k y)
 =\partial_i^\alpha x\circ_{k-1}\partial_i^\alpha y
 \ \text{for $i<k$},\\
 &\partial_k^-(x\circ_k y)=\partial_k^- x,\\
 &\partial_k^+(x\circ_k y)=\partial_k^+ y,\\
 &\partial_i^\alpha(x\circ_k y)
 =\partial_i^\alpha x\circ_k\partial_i^\alpha y
 \ \text{for $i>k$}.
\end{align*}
\end{proposition}

\begin{proof}
This follows straightforwardly from the definition: if $i<k$ then
$\partial_i^\alpha(x\circ_k y)=\partial_i^\alpha
x\circ_{k-1}\partial_i^\alpha y$ because
$$\partial_i^\alpha\partial_k^- G_k(x,y)
 =\partial_{k-1}^-\partial_i^\alpha G_k(x,y)
 =\partial_{k-1}^-
 G_{k-1}(\partial_i^\alpha x,\partial_i^\alpha y),$$
etc.
\end{proof}

\begin{proposition} \label{4I}
The composites have the properties
\begin{align*}
 &\epsilon_k\partial_k^- x\circ_k x
 =x
 =x\circ_k\epsilon_k\partial_k^+ x,\\
 &\Gamma_k^+ x\circ_k\Gamma_k^- x=\epsilon_{k+1}x,\\
 &\Gamma_k^+ x\circ_{k+1}\Gamma_k^- x=\epsilon_k x.
\end{align*}
\end{proposition}

\begin{proof}
There are composers $G_k(\epsilon_k\partial_k^- x,x)$ and
$G_k(x,\epsilon_k\partial_k^+ x)$ because
$$\partial_k^+\epsilon_k\partial_k^- x=\partial_k^- x,\qquad
 \partial_k^+ x=\partial_k^-\epsilon_k\partial_k^+ x.$$
An inductive argument shows that $\partial_i^\alpha
G_k(\epsilon_k\partial_k^- x,x)=\partial_i^\alpha\epsilon_k x$ for
$\partial_i^\alpha\neq\partial_k^-$, so that $\epsilon_k x$ is a
thin filler for the admissible box whose unique thin filler is
$G_k(\epsilon_k\partial_k^- x,x)$. Therefore
$G_k(\epsilon_k\partial_k^- x,x)=\epsilon_k x$ (compare Figures
\ref{F6} and~\ref{F8}). By a similar argument,
$G_k(x,\epsilon_k\partial_k^+ x)=\Gamma_k^- x$ (compare Figures
\ref{F7} and~\ref{F8}). Applying~$\partial_k^-$ now gives
$\epsilon_k\partial_k^- x\circ_k x=x$ and
$x\circ_k\epsilon_k\partial_k^+ x=x$.

It is clear that $\Gamma_k^+ x\circ_k\Gamma_k^- x$ exists. From
Proposition~\ref{4E} it is thin, and an inductive argument shows
that it is a filler for the shell whose unique thin filler is
$\epsilon_{k+1}x$. Therefore $\Gamma_k^+ x\circ_k\Gamma_k^-
x=\epsilon_{k+1}x$. Similarly $\Gamma_k^+ x\circ_{k+1}\Gamma_k^-
x=\epsilon_k x$.
\end{proof}

\begin{proposition} \label{4J}
If $x\circ_k y$ is defined then
\begin{align*}
 &\epsilon_j(x\circ_k y)=\epsilon_j x\circ_{k+1}\epsilon_j y
 \ \text{for $j\leq k$},\\
 &\epsilon_j(x\circ_k y)=\epsilon_j x\circ_k\epsilon_j y
 \ \text{for $j>k$},\\
 &\Gamma_j^\beta(x\circ_k y)=\Gamma_j^\beta x\circ_{k+1}\Gamma_j^\beta y
 \ \text{for $j<k$},\\
 &\Gamma_k^-(x\circ_k y)
 =(\Gamma_k^- x\circ_k\epsilon_{k+1}y)\circ_{k+1}\Gamma_k^- y
 =(\Gamma_k^- x\circ_{k+1}\epsilon_k y)\circ_k\Gamma_k^- y,\\
 &\Gamma_k^+(x\circ_k y)
 =\Gamma_k^+ x\circ_k(\epsilon_k x\circ_{k+1}\Gamma_k^+ y)
 =\Gamma_k^+ x\circ_{k+1}(\epsilon_{k+1}x\circ_k\Gamma_k^+ y),\\
 &\Gamma_j^\beta(x\circ_k y)=\Gamma_j^\beta x\circ_k\Gamma_j^\beta y
 \ \text{for $j>k$}.
\end{align*}
\end{proposition}

\begin{proof}
In each equality the first expression is defined as the unique
thin filler of some shell and the other expressions are
well-defined and thin. It therefore suffices to show that the
other expressions are also fillers for the appropriate shells, and
this is done by inductive arguments.
\end{proof}

\begin{proposition} \label{4K}
If $k\neq l$, then
$$(x\circ_k y)\circ_l(z\circ_k w)
 =(x\circ_l z)\circ_k(y\circ_l w)$$
whenever both sides are defined.
\end{proposition}

\begin{proof}
For definiteness, suppose that $k>l$. By an inductive argument
one shows that
$$G_k(x,y)\circ_l G_k(z,w)=G_k(x\circ_l z,y\circ_l w);$$
indeed the composite on the left exists and is a thin filler for
the box opposite~$\partial_k^-$ whose unique thin filler is the
expression on the right. The result then follows by
applying~$\partial_k^-$ to both sides.
\end{proof}

\begin{proposition} \label{4L}
If $x,y,z\in X_n$ are such that $\partial_k^+ x=\partial_k^- y$
and $\partial_k^+ y=\partial_k^- z$, then
$$(x\circ_k y)\circ_k z=x\circ_k(y\circ_k z).$$
\end{proposition}

\begin{figure}
 $$\xymatrix{
 \ar@{-}[]+0;[rrrr]+0^{\epsilon_1\partial_1^+ z}&&&&\\
 &&*-=0{G_1(x\circ_1 y,z)}&*i{G_1(x,y\circ_1 z)}\\
 \ar@{.}+0;[uu]+0^{(x\circ_1 y)\circ_1 z}
 \ar@{-}[]+0;[rrrr]+0_{x\circ_1 y}&&&&
 \ar@{-}[]+0;[uu]+0_{z}
 \\
 \ar@{-}[]+0;[rr]+0^{\epsilon_1\partial_1^+ z}&&
 \ar@{-}[]+0;[rr]+0^{\epsilon_1\partial_1^+ z}&&\\
 &{G_1(x,y\circ_1 z)}&&{G_1(y,z)}\\
 \ar@{.}[]+0;[uu]+0^{x\circ_1(y\circ_1 z)}
 \ar@{-}[]+0;[rr]+0_x&&
 \ar@{-}[]+0;[uu]+0_{y\circ_1 z}
 \ar@{-}[]+0;[rr]+0_{y}&&
 \ar@{-}[]+0;[uu]+0_{z}
 }$$
 \caption{The equality
 $G_1(x\circ_1 y,z)=G_1(x,y)\circ_1 G_1(x,y\circ_1 z)$}
 \label{F9}
\end{figure}
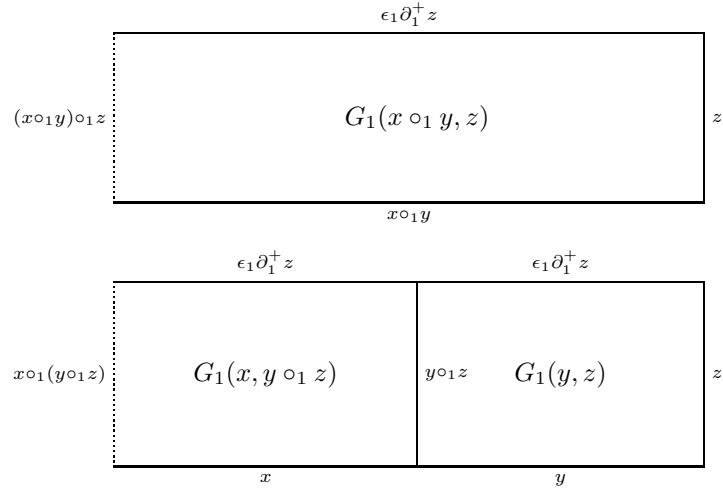

\begin{proof}
We first show that
$$G_k(x\circ_k y,z)=G_k(x,y\circ_k z)\circ_k G_k(y,z)$$
by the usual inductive argument: the expression on the right
exists and is a thin filler for the admissible box
opposite~$\partial_k^-$ whose unique thin filler is the expression
on the left. The case $k=n=1$ is shown in Figure~\ref{F9}. We then
get $\partial_k^-G_k(x\circ_k y,z)=\partial_k^-G_k(x,y\circ_k z)$,
which means that $(x\circ_k y)\circ_k z=x\circ_k(y\circ_k z)$ as
required.
\end{proof}

We have now verified all the axioms of \cite{ABS}, so we have
proved the following result.

\begin{theorem} \label{4M}
If $X$~is a complete stratified precubical set then the induced
degeneracies, connections and compositions make~$X$ into the
cubical nerve of an $\omega$-category.
\end{theorem}

\section{The equivalence} \label{S5}

We have shown that an $\omega$-category nerve structure on a
precubical set induces a complete stratification
(Theorem~\ref{3N}) and that a complete stratification induces an
$\omega$-category nerve structure (Theorem~\ref{4M}). In this
section we complete the proof of Theorem~\ref{2N} by showing that
the two processes are mutually inverse.

We begin by recalling some properties of nerves from~\cite{ABS}.

\begin{proposition} \label{5A}
Let $X$ be the cubical nerve of an $\omega$-category~$C$, so that
$X$~has degeneracies, connections and compositions and an induced
complete stratification. Then degeneracies and connections are
thin, and composites of thin elements are thin.
\end{proposition}

\begin{proof}
Let $C(n)$ be the sub-$\omega$-category of~$C$ consisting of the
elements which are identities for~$\comp_n$ (see
Proposition~\ref{3A}), and let $X(n)$ be the nerve of $C(n)$; then
$X(n)$ is a sub-precubical set of~$X$ closed under degeneracies,
connections and compositions. We have $X_n=X(n)_n$, and the thin
elements of~$X_n$ are precisely the members of $X(n-1)_n$. If an
$n$-cube~$x$ is a degeneracy or connection of an $(n-1)$-cube~$y$,
then $y\in X_{n-1}=X(n-1)_{n-1}$, so $x\in X(n-1)_n$ and $x$~is
therefore thin. If an $n$-cube~$x$ is a composite of thin
$n$-cubes, then $x$~is thin because $X(n-1)$ is closed under
composition. This completes the proof.
\end{proof}

\begin{proposition} \label{5B}
Let $X$ be the cubical nerve of an $\omega$-category. Then there
are operations
$$\psi_1,\ldots,\psi_{n-1}\colon X_n\to X_n$$
such that
\begin{align*}
 &\psi_k x
 =\Gamma_k^+\partial_{k+1}^-x
 \circ_{k+1}x
 \circ_{k+1}\Gamma_k^-\partial_{k+1}^+x,\\
 &\partial_{k+1}^\alpha\psi_k x
 =\epsilon_k\partial_k^\alpha\partial_{k+1}^\alpha x,\\
 &\partial_i^\alpha\psi_k x=\psi_k\partial_i^\alpha x
 \ \text{for $i>k+1$},\\
 &\psi_k\epsilon_{k+1}x=\epsilon_k x,\\
 &x
 =(\epsilon_k\partial_k^-x
 \circ_{k+1}\Gamma_k^+\partial_{k+1}^+x)
 \circ_k\psi_k x
 \circ_k(\Gamma_k^-\partial_{k+1}^-x
 \circ_{k+1}\epsilon_k\partial_k^+x).
\end{align*}
\end{proposition}

\begin{proof}
We know from~\cite{ABS} that $X$~is a precubical set with
operations satisfying the conditions of Definitions
\ref{4A}--\ref{4C}, Propositions \ref{4H}--\ref{4I} and
Propositions \ref{4K}--\ref{4L}. Using these conditions, it is
easy to check that the composite on the right side of the first
equality exists, and one can therefore use this equality to define
$\psi_k x$. It is then straightforward to verify the next three
equalities. As to the last equality, we have
\begin{align*}
 x
 &=x\circ_{k+1}\epsilon_{k+1}\partial_{k+1}^+x\\
 &=(\epsilon_k\partial_k^-x\circ_k x)
 \circ_{k+1}
 (\Gamma_k^+\partial_{k+1}^+x\circ_k\Gamma_k^-\partial_{k+1}^+x)\\
 &=(\epsilon_k\partial_k^-x\circ_{k+1}\Gamma_k^+\partial_{k+1}^+x)
 \circ_k
 (x\circ_{k+1}\Gamma_k^-\partial_{k+1}^+x)
\end{align*}
and
\begin{align*}
 x\circ_{k+1}\Gamma_k^-\partial_{k+1}^+x
 &=\epsilon_{k+1}\partial_{k+1}^-x\circ_{k+1}
 (x\circ_{k+1}\Gamma_k^-\partial_{k+1}^+x)\\
 &=(\Gamma_k^+\partial_{k+1}^-x\circ_k\Gamma_k^-\partial_{k+1}^-x)
 \circ_{k+1}
 [(x\circ_{k+1}\Gamma_k^-\partial_{k+1}^+x)
 \circ_k\epsilon_k\partial_k^+x]\\
 &=[\Gamma_k^+\partial_{k+1}^-x
 \circ_{k+1}(x\circ_{k+1}\Gamma_k^-\partial_{k+1}^+x)]
 \circ_k
 (\Gamma_k^-\partial_{k+1}^-x\circ_{k+1}\epsilon_k\partial_k^+x)\\
 &=\psi_k x
 \circ_k(\Gamma_k^-\partial_{k+1}^-x\circ_{k+1}\epsilon_k\partial_k^+x),
\end{align*}
from which the result follows.
\end{proof}

\begin{proposition} \label{5C}
Let $x$ be an $n$-cube in the cubical nerve of an
$\omega$-category with $n>0$. Then there is an $n$-cube $\Psi x$
such that $\Psi x$ can be obtained by composing~$x$ with
connections, such that $x$~can be obtained by composing $\Psi x$
with degeneracies and connections, and such that
$$\partial_i^\alpha\Psi x
 =\epsilon_1\partial_{i-1}^\alpha\partial_1^+\Psi x
 \ \text{for $i>1$}.$$
\end{proposition}

\begin{proof}
Let
$$\Psi x=\psi_1\psi_2\ldots\psi_{n-1}x,$$
where the~$\psi_k$ are as in Proposition~\ref{5B}. From the first
equality in Proposition~\ref{5B}, $\Psi x$ is a composite of~$x$
with connections; from the last equality in Proposition~\ref{5B},
$x$~is a composite of $\Psi x$ with degeneracies and connections.
For $i>1$, the middle equalities in Proposition~\ref{5B} give
\begin{align*}
 \partial_i^\alpha\Psi x
 &=\partial_i^\alpha(\psi_1\ldots\psi_{i-2})\psi_{i-1}
  (\psi_i\ldots\psi_{n-1})x\\
 &=(\psi_1\ldots\psi_{i-2})\partial_i^\alpha\psi_{i-1}
  (\psi_i\ldots\psi_{n-1})x\\
 &=(\psi_1\ldots\psi_{i-2})
 \epsilon_{i-1}\partial_{i-1}^\alpha\partial_i^\alpha
  (\psi_i\ldots\psi_{n-1})x\\
 &=\epsilon_1\partial_{i-1}^\alpha\partial_i^\alpha
  (\psi_i\ldots\psi_{n-1})x\\
 &=\epsilon_1 y,
\end{align*}
say, we then get
$$y
 =\partial_1^+\epsilon_1 y
 =\partial_1^+\partial_i^\alpha\Psi x
 =\partial_{i-1}^\alpha\partial_1^+\Psi x,$$
and we deduce that $\partial_i^\alpha\Psi x=\epsilon_1
y=\epsilon_1\partial_{i-1}^\alpha\partial_1^+\Psi x$. This
completes the proof.
\end{proof}

We now give the two results showing that the functors are mutually
inverse.

\begin{theorem} \label{5D}
Let $X$ be a complete stratified precubical set. Then the
stratification on~$X$ obtained from its structure as the cubical
nerve of an $\omega$-category is the same as the original
stratification.
\end{theorem}

\begin{proof}
We use the method of Higgins \cite{Higg}. Let $x$ be an $n$-cube
in~$X$ with $n>0$, and let $\Psi x$ be as in Proposition~\ref{5C}.
From Propositions \ref{5A} and \ref{5C} we see that in the
$\omega$-category stratification $x$~is thin if and only if $\Psi
x$ is thin. But Proposition~\ref{5A} is also true in the original
stratification: degeneracies and connections are thin by
construction, and composites of thin elements are thin by
Proposition~\ref{4E}. Hence, in the original stratification, it
also follows from Proposition~\ref{5C} that $x$~is thin if and
only if $\Psi x$ is thin. It therefore suffices to show that $\Psi
x$ is thin in the $\omega$-category stratification if and only if
it is thin in the original stratification.

From Proposition~\ref{5C} we see that $\partial_i^\alpha\Psi
x=\partial_i^\alpha\epsilon_1\partial_1^+\Psi x$ for $i>1$, and we
also have $\partial_1^+\Psi
x=\partial_1^+\epsilon_1\partial_1^+\Psi x$, so $\Psi x$ and
$\epsilon_1\partial_1^+\Psi x$ are fillers for the same
$n$-box~$b$ opposite~$\partial_1^-$. In each of the
stratifications degeneracies are thin, so $b$~is admissible as in
Example~\ref{2J} (see the justification of Definition~\ref{4A}),
and $b$~therefore has a unique thin filler. Since degeneracies are
thin in each stratification, the unique thin filler of~$b$ in each
stratification is given by $\epsilon_1\partial_1^+\Psi x$. In each
stratification, it follows that $\Psi x$ is thin if and only if
$\Psi x=\epsilon_1\partial_1^+\Psi x$. Therefore $\Psi x$ is thin
in the $\omega$-category stratification if and only if it is thin
in the original stratification. This completes the proof.
\end{proof}

\begin{theorem} \label{5E}
Let $X$ be the cubical nerve of an $\omega$-category. Then the
cubical nerve structure obtained from the induced stratification
is the same as the original cubical nerve structure.
\end{theorem}

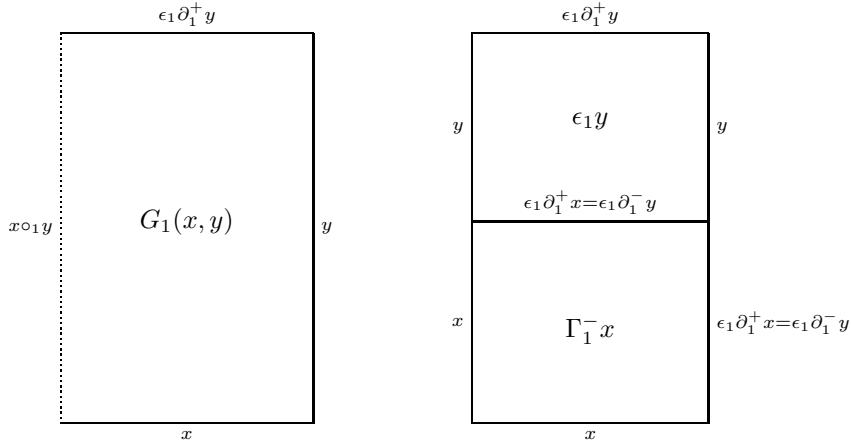
\begin{figure}
 $$\xymatrix{
 \ar@{-}[]+0;[rr]+0^{\epsilon_1\partial_1^+ y}
 &&&&\ar@{-}[]+0;[rr]+0^{\epsilon_1\partial_1^+ y}&&\\
 &&&&&\epsilon_1 y\\
 &G_1(x,y)&&&\ar@{-}[]+0;[uu]+0^{y}
 \ar@{-}[]+0;[rr]+0^{\epsilon_1\partial_1^+ x=\epsilon_1\partial_1^- y}&
 *i{G_1(x,y)}&
 \ar@{-}[]+0;[uu]+0_{y}\\
 &&&&&\Gamma_1^- x\\
 \ar@{.}[]+0;[uuuu]+0^{x\circ_1 y}
 \ar@{-}[]+0;[rr]+0_{x}&&
 \ar@{-}[]+0;[uuuu]+0_{y}&&
 \ar@{-}[]+0;[uu]+0^{x}
 \ar@{-}[]+0;[rr]+0_{x}&&
 \ar@{-}[]+0;[uu]+0_{\epsilon_1\partial_1^+ x=\epsilon_1\partial_1^- y}
 }$$
 \caption{The equality $G_1(x,y)=\Gamma_1^- x\circ_2\epsilon_1 y$}
 \label{F10}
\end{figure}

\begin{proof}
We must show that the degeneracies, connections and compositions
constructed from the stratification are the same as the original
degeneracies, connections and compositions. Now the original
degeneracies are thin by Proposition~\ref{5A}, so they satisfy the
conditions of Definition~\ref{4A}, and it follows that they are
the same as the degeneracies constructed from the stratification.
The same argument applies to connections. As to compositions, let
$x$~and~$y$ be $n$-cubes such that $\partial_k^+ x=\partial_k^- y$
for some~$k$. In the original structure one can check that there
is a composite $\Gamma_k^- x\circ_{k+1}\epsilon_k y$, and this
composite is thin by Proposition~\ref{5A}. By an inductive
argument one finds that this composite satisfies the conditions
for $G_k(x,y)$ in Definition~\ref{4C}; see Figure~\ref{F10} for
the case $n=k=1$. One therefore gets
$$G_k(x,y)=\Gamma_k^- x\circ_{k+1}\epsilon_k y,$$
where the left side is the composer constructed from the
stratification and the right side is the composite in the original
structure. Applying~$\partial_k^-$ to both sides now shows that
the composite $x\circ_k y$ constructed form the stratification is
the same as the original composite $\partial_k^-\Gamma_k^-
x\circ_k\partial_k^-\epsilon_k y=x\circ_k y$. This completes the
proof.
\end{proof}

\section{Proof of Theorem 3.14} \label{S6}

Let $\partial_k^\gamma$ be a face operation on $n$-cubes; we must
construct a factorisation of $d_{n-1}^{(-)^{k-1}\gamma}\langle
u_n\rangle$ in $\nu I^n$ with certain properties. We will write
$$\sigma
 =\check\partial_k^\gamma u_{n-1}
 =u_{k-1}\otimes\partial^\gamma u_1\otimes u_{n-k},$$
so that an $n$-box opposite~$\partial_k^\gamma$ is a morphism on
$B(\sigma)$.

Consider the self-map of a geometric $n$-cube got by projection
onto the face corresponding to~$\sigma$. It is a cellular map
cellularly homotopic to the identity, so it induces a chain
endomorphism of~$I^n$ chain homotopic to the identity. To be
explicit, we get the following two results.

\begin{proposition} \label{6Y}
There is a chain map $f\colon I^n\to I^n$ given for $x$ a chain
in~$I^{k-1}$ and for $y$ a chain in~$I^{n-k}$ by
\begin{align*}
 &f(x\otimes u_1\otimes y)=0,\\
 &f(x\otimes\partial^\gamma u_1\otimes y)
 =f(x\otimes\partial^{-\gamma}u_1\otimes y)
 =x\otimes\partial^\gamma u_1\otimes y,
\end{align*}
such that
$$f\langle u_n\rangle_i^\alpha=\langle\sigma\rangle_i^\alpha.$$
\end{proposition}

\begin{proof}
It is straightforward to check that there is a chain map~$f$ as
given, and it follows from Proposition~\ref{3F} that $f\langle
u_n\rangle_i^\alpha=\langle\sigma\rangle_i^\alpha$.
\end{proof}

\begin{proposition} \label{6Z}
There are abelian group homomorphisms $D\colon I^n\to I^n$ of
degree~$1$, given for $x$ an $i$-chain in~$I^{k-1}$ and for $y$ a
chain in~$I^{n-k}$ by
\begin{align*}
 &D(x\otimes u_1\otimes y)=0,\\
 &D(x\otimes\partial^\gamma u_1\otimes y)=0,\\
 &D(x\otimes\partial^{-\gamma}u_1\otimes y)
 =-(-)^i\gamma(x\otimes u_1\otimes y),
\end{align*}
such that
$$\partial D+D\partial=\id-f.$$
\end{proposition}

\begin{proof}
This is a straightforward computation.
\end{proof}

The chain homotopy~$D$ is related to the precubical operations
complementary to~$\partial_k^\gamma$ as follows.

\begin{proposition} \label{6B}
The chains $D\langle u_n\rangle_{q-1}^+$ and $-D\langle
u_n\rangle_{q-1}^-$ are sums of basis elements
$$\check\partial_{i(n-q)}^{\alpha(n-q)}\ldots
 \check\partial_{i(1)}^{\alpha(1)} u_q$$
corresponding to precubical operations
$\partial_{i(1)}^{\alpha(1)}\ldots\partial_{i(n-q)}^{\alpha(n-q)}$
complementary to~$\partial_k^\gamma$.
\end{proposition}

\begin{proof}
This follows from Proposition~\ref{3F}, according to which
$\langle u_n\rangle_{q-1}^\alpha$ is the sum of the basis elements
$$\check\partial_{i(n-q+1)}^{\alpha(n-q+1)}\ldots
 \check\partial_{i(1)}^{\alpha(1)} u_q$$
such that $i(1)<i(2)<\ldots<i(n-q+1)$ and such that
$(-)^{i(r)-r}\alpha(r)=\alpha$ for all~$r$. Applying~$D$ picks out
the terms involving~$\check\partial_k^{-\gamma}$, omits the
factor~$\check\partial_k^{-\gamma}$, and multiplies by~$\alpha$.
This makes $D\langle u_n\rangle_{q-1}^+$ and $-D\langle
u_n\rangle_{q-1}^-$ into sums of basis elements corresponding to
precubical operations complementary to~$\partial_k^\gamma$, as
required.
\end{proof}

The factors~$A_q^\beta$ of $d_{n-1}^{(-)^{k-1}\gamma}\langle
u_n\rangle$ are defined as follows.

\begin{proposition} \label{6X}
There are elements $A_q^-$~and~$A_q^+$ of $\nu I^n$ for $1\leq
q\leq n-1$ given as double sequences by the formulae
\begin{align*}
 &(A_q^\beta)_i^\alpha=\langle u_n\rangle_i^\alpha
 \ \text{for $i<q-1$},\\
 &(A_q^\beta)_{q-1}^\beta=\langle u_n\rangle_{q-1}^\beta,\\
 &(A_q^\beta)_{q-1}^{-\beta}
 =(\id-\partial D)\langle u_n\rangle_{q-1}^\beta,\\
 &(A_q^\beta)_q^\alpha=\beta D\langle u_n\rangle_{q-1}^\beta,\\
 &(A_q^\beta)_i^\alpha=0\ \text{for $i>q$}.
\end{align*}
\end{proposition}

\begin{proof}
We must verify the conditions of Definition~\ref{3D}. Note that
$(A_q^\beta)_{q-1}^{-\beta}$~can be written in the form
$$(f+D\partial)\langle u_n\rangle_{q-1}^\beta
 =\langle\sigma\rangle_{q-1}^\beta
 +D\langle u_n\rangle_{q-2}^+
 -D\langle u_n\rangle_{q-2}^-$$
(interpret $\langle u_n\rangle_{-1}^\alpha$ as zero in the case
$q=1$). Using Proposition~\ref{6B} where necessary, we see that
$(A_q^\beta)_i^\alpha$ is a sum of $i$-dimensional basis elements
in~$I^n$. From the expressions in the statement of the
proposition, it is easy to check that
$\epsilon(A_q^\beta)_0^\alpha=1$ and $(A_q^\beta)_i^+
-(A_q^\beta)_i^- =\partial(A_q^\beta)_{i+1}^\alpha$. This
completes the proof.
\end{proof}

We must now show that the elements~$A_q^\beta$ have the properties
stated in Theorem~\ref{3O}. We begin with the following
observation.

\begin{proposition} \label{6W}
The elements~$A_q^\beta$ are members of $\nu B(\sigma)$ such that
$b(A_q^\beta)$ is an identity for~$\comp_{q-1}$ whenever $b$~is an
admissible $n$-box opposite~$\partial_k^\gamma$.
\end{proposition}

\begin{proof}
The~$A_q^\beta$ are members of $\nu B(\sigma)$ because the basis
element~$\sigma$ is never a term in~$(A_q^\beta)_{n-1}^\alpha$ and
because $u_n$~is never a term in~$(A_q^\beta)_n^\alpha$.

Now let $b$ be an admissible $n$-box opposite~$\partial_k^\gamma$.
By Proposition~\ref{6B} and the definition of admissible box,
$b\langle\tau\rangle$ is an identity for~$\comp_{q-1}$ whenever
$\tau$~is a term in $\beta D\langle u_n\rangle_{q-1}^\beta$. But
$(A_q^\beta)_i^\alpha=0$ for $i>q$ and
$$(A_q^\beta)_q^-
 =(A_q^\beta)_q^+
 =\beta D\langle u_n\rangle_{q-1}^\beta;$$
hence, by Theorem~\ref{3I}, $A_q^\beta$~is a composite of
atoms~$\langle\tau\rangle$ such that $\tau$~is a term in $\beta
D\langle u_n\rangle_{q-1}^\beta$ or $\tau$~has dimension less
than~$q$. It follows that $b(A_q^\beta)$ is a composite of
identities for~$\comp_{q-1}$. By Proposition~\ref{3A},
$b(A_q^\beta)$ itself is an identity for~$\comp_{q-1}$.

This completes the proof.
\end{proof}

Next we show how to compose the elements~$A_q^\beta$.

\begin{proposition} \label{6V}
There are elements~$A_q$ in $\nu I^n$ for $0\leq q\leq n-1$ given
inductively by $A_0=\langle\sigma\rangle$ and by
$$A_q=A_q^-\comp_{q-1} A_{q-1}\comp A_q^+\ \text{for $1\leq q\leq n-1$}$$
such that
\begin{align*}
 &(A_q)_i^\alpha=\langle u_n\rangle_i^\alpha
 \ \text{for $i<q$},\\
 &(A_q)_q^\alpha=(\id-\partial D)\langle u_n\rangle_q^\alpha,\\
 &(A_q)_i^\alpha=\langle \sigma\rangle_i^\alpha
 \ \text{for $i>q$}.
\end{align*}
\end{proposition}

\begin{proof}
The proof is by induction. To begin, let
$A_0=\langle\sigma\rangle$; then $A_0$~is certainly a member of
$\nu I^n$ such that the $(A_0)_i^\alpha$ are as described, since
$$(\id-\partial D)\langle u_n\rangle_0^\alpha
 =(f+D\partial)\langle u_n\rangle_0^\alpha
 =f\langle u_n\rangle_0^\alpha
 =\langle\sigma\rangle_0^\alpha.$$

For the inductive step, suppose that there is a member~$A_{q-1}$
of $\nu I^n$ as described. Then one finds that $d_{q-1}^+
A_q^-=d_{q-1}^- A_{q-1}$ and $d_{q-1}^+ A_{q-1}=d_{q-1}^- A_q^+$,
so one can define~$A_q$ in $\nu I^n$ as the composite
$$A_q
 =A_q^-\comp_{q-1} A_{q-1}\comp A_q^+
 =A_q^- -d_{q-1}^- A_{q-1}+A_{q-1}-d_{q-1}^+ A_{q-1}+A_q^+.$$
It now follows from the inductive hypothesis that
the~$(A_q)_i^\alpha$ are as required; for the case $i=q$ note that
$$-D\langle u_n\rangle_{q-1}^-
+\langle\sigma\rangle_q^\alpha
 +D\langle u_n\rangle_{q-1}^+
 =(f+D\partial)\langle u_n\rangle_q^\alpha
 =(\id-\partial D)\langle u_n\rangle_q^\alpha.$$

This completes the proof.
\end{proof}

Finally, we complete the proof of Theorem~\ref{3O} by proving the
following result.

\begin{proposition} \label{6U}
There is an equality
$$A_{n-1}=d_{n-1}^{(-)^{k-1}\gamma}\langle u_n\rangle.$$
\end{proposition}

\begin{proof}
We have $(A_{n-1})_i^\alpha=\langle u_n\rangle_i^\alpha$ for
$i<n-1$ and $(A_{n-1})_i^\alpha=\langle\sigma\rangle_i^\alpha=0$
for $i>n-1$, from which it follows that
$$(A_{n-1})_{n-1}^+ -(A_{n-1})_{n-1}^-
 =\partial(A_{n-1})_n^-
 =0.$$
It therefore suffices to show that
$(A_{n-1})_{n-1}^{(-)^{k-1}\gamma}=\langle
u_n\rangle_{n-1}^{(-)^{k-1}\gamma}$. But
$$(A_{n-1})_{n-1}^{(-)^{k-1}\gamma}
 =(\id-\partial D)\langle u_n\rangle_{n-1}^{(-)^{k-1}\gamma}
 =\langle u_n\rangle_{n-1}^{(-)^{k-1}\gamma}$$
because $u_{k-1}\otimes\partial^{-\gamma}u_1\otimes u_{n-k}$ is
not a term in $\langle u_n\rangle_{n-1}^{(-)^{k-1}\gamma}$. This
completes the proof.
\end{proof}

\section{A thin filler with a single non-thin face} \label{S7}

In this section we exhibit a thin $3$-cube~$x$ in the nerve of an
$\omega$-category such that $x$~has exactly one non-thin $2$-face;
moreover, $x$~is the thin filler of an admissible box. It complies
with the final condition of Definition~\ref{2M}, because the
non-thin face is $\partial_1^- x$ and the box is opposite a
different operation~$\partial_2^-$, but it shows that this
condition in Definition~\ref{2M} cannot be weakened.

We need an $\omega$-category containing elements with certain
properties.

\begin{proposition} \label{7A}
There is an $\omega$-category with elements $A$~and~$b$ such that
$A$~is an identity for~$\comp_2$ but not for~$\comp_1$, such that
$b$~is an identity for~$\comp_1$, and such that $A\comp_0 b$
exists and is an identity for~$\comp_1$.
\end{proposition}

\begin{proof}
We take the $\omega$-category to be the $2$-category of small
categories; thus the identities for~$\comp_0$ are the categories,
the identities for~$\comp_1$ are the functors, the identities
for~$\comp_2$ are the natural transformations, and every element
is an identity for~$\comp_3$. Let $A$ be a non-identity natural
transformation and $b$ be a functor such that $A\comp_0 b$ exists
and is an identity natural transformation. Then $A$~and~$b$ have
the required properties.
\end{proof}

\begin{theorem} \label{7B}
There is an $\omega$-category with a $3$-cube~$x$ in its cubical
nerve such that $x$~is the thin filler of an admissible box
opposite~$\partial_2^-$, such that $\partial_1^- x$ is not thin,
and such that $\partial_i^\alpha x$ is thin for
$\partial_i^\alpha\neq\partial_1^-$.
\end{theorem}

\begin{figure}
 $$\xymatrix@!C{
 &&\ar[rrrr]^{a^-}
 &&&&\ar[ddrr]^{b}
 \\
 \\
 \ar[uurr]^{=}
 \ar@{}[uurrrrrr]|{A\Downarrow}
 \ar[rrrr]^{a^+}
 \ar@{}[ddrrrrrr]|{a^+\Downarrow}
 \ar[ddrr]_{=}
 &&&&\ar[uurr]_{=}
 \ar@{}[rrrr]|{b\Downarrow}
 \ar[ddrr]^{=}
 &&&&\\
 \\
 &&\ar[rrrr]_{a^+}
 &&&&\ar[uurr]_{b}
 \\
 &&&&\\
 &&\ar[rrrr]^{a^-}
 \ar@{}[ddrrrrrr]|{a^-\comp_0 b\Downarrow}
 \ar[ddrr]_{a^-\comp_0 b}
 &&&&\ar[ddrr]^{b}
 \\
 \\
 \ar[uurr]^{=}
 \ar@{}[rrrr]|{A\comp_0 b\Downarrow}
 \ar[ddrr]_{=}
 &&&&\ar[rrrr]^{=}
 &&&&\\
 \\
 &&\ar[uurr]^{a^+\comp_0 b}
 \ar@{}[uurrrrrr]|{a^+\comp_0 b\Downarrow}
 \ar[rrrr]_{a^+}
 &&&&\ar[uurr]_{b}
 }$$
 \caption{The shell $s$ in the proof of Theorem \ref{7B}}
 \label{F11}
\end{figure}
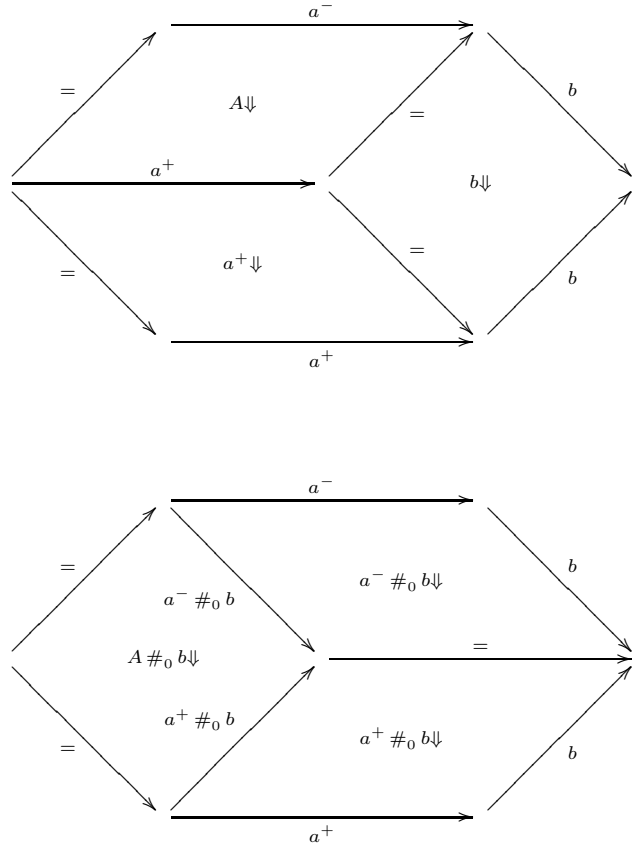

\begin{proof}
By Proposition~\ref{3L}, to construct the thin $3$-cube~$x$ it
suffices to construct a $3$-shell~$s$ such that $s(d_2^-\langle
u_3\rangle)=s(d_2^+\langle u_3\rangle)$. By Proposition~\ref{3H},
to construct the shell~$s$, it suffices to assign
$\omega$-category elements to the atoms of dimension less than~$3$
in~$I^3$ such that the restriction to each $2$-face is as shown in
Figure~\ref{F1}. Take an $\omega$-category with elements
$A$~and~$b$ as in Proposition~\ref{7A}, let $a^-=d_1^- A$, and let
$a^+=d_1^+ A$. One can then check that there is a shell~$s$ as
shown in Figure~\ref{F11}, where the equality signs denote
identities for~$\comp_0$ and where the $2$-faces are positioned as
in Figure~\ref{F4}; that is,
\begin{align*}
 &x\langle\check\partial_1^- u_2\rangle=A,\
 x\langle\check\partial_2^+ u_2\rangle=b,\
 x\langle\check\partial_3^- u_2\rangle=a^+,\\
 &x\langle\check\partial_3^+ u_2\rangle=a^-\comp_0 b,\
 x\langle\check\partial_2^- u_2\rangle=A\comp_0 b,\
 x\langle\check\partial_1^+ u_2\rangle=a^+\comp_0 b.
\end{align*}
From the formulae for $d_2^-\langle u_3\rangle$ and $d_2^+\langle
u_3\rangle$ in Example~\ref{3Y}, we see that $s(d_2^-\langle
u_3\rangle)$ and $s(d_2^+\langle u_3\rangle)$ are both equal to
$A\comp_0 b$, so $s$~has a thin filler~$x$. Since $a^-$, $a^+$,
$b$ and $A\comp_0 b$ are identities for~$\comp_1$ and since $A$~is
not an identity for~$\comp_1$, it follows that $\partial_1^- x$ is
the unique non-thin $2$-face of~$x$. In particular $\partial_1^+
x$ and $\partial_3^+ x$ are thin; so also is
$\partial_1^+\partial_3^+ x$, which is given by the common edge of
$\partial_1^+ x$ and $\partial_3^+ x$. This means that
$x$~restricts to an admissible box opposite~$\partial_2^-$, so
$x$~is the thin filler of an admissible box
opposite~$\partial_2^-$. This completes the proof.
\end{proof}

\providecommand{\bysame}{\leavevmode\hbox
to3em{\hrulefill}\thinspace}
\providecommand{\MR}{\relax\ifhmode\unskip\space\fi MR }
\providecommand{\MRhref}[2]{%
  \href{http://www.ams.org/mathscinet-getitem?mr=#1}{#2}
} \providecommand{\href}[2]{#2}

\end{document}